\documentclass[10pt]{article}

\usepackage{amsfonts,amsmath,amssymb,amsthm}
\usepackage{doublespace,fullpage}
\usepackage{graphics,times,kmac}
\usepackage{epic,eepic}
\usepackage{calc}

\usepackage[dvipdfm]{hyperref}

\newcommand{\mapone}{$1.996\pm 0.0071$}
\newcommand{\maponerange}{[10^{-3.5},10^{-1}]}

\newcommand{\maptwo}{$0.97\pm 0.033$}
\newcommand{\maptworange}{[10^{-3},10^{-1}]}

\newcommand{\mapthreea}{$0.52\pm 0.056$}
\newcommand{\mapthreearange}{[10^{-5},10^{-2.5}]}
\newcommand{\mapthreeb}{$0.46\pm 0.050$}
\newcommand{\mapthreebrange}{[10^{-5},10^{-3}]}

\newcommand{\mapfoura}{$0.31\pm 0.028$}
\newcommand{\mapfourarange}{[10^{-4},10^{-1.5}]}
\newcommand{\mapfourb}{$0.53\pm 0.056$}
\newcommand{\mapfourbrange}{[10^{-4},10^{-2}]}
\newcommand{\mapfourc}{$0.17\pm 0.033$}
\newcommand{\mapfourcrange}{[10^{-4},10^{-1.5}]}

\newcommand{\mapfive}{$1.40\pm 0.022$}
\newcommand{\mapfiverange}{[10^{-2.5},10^{-1}]}

\newcommand{\putgraph}[1]{\includegraphics[0in,0in][4in,3in]{pix/#1}}
\newcommand{\onenorm}[1]{\norm{#1}_1}

\newcommand{\tcorr}{\varsub{\tau}{corr}}
\newcommand{\musrb}{\varsub{\mu}{SRB}}

\newcounter{hour}
\newcounter{minute}
\newcounter{timenow}
\newcommand{\now}{
\ifnum \value{minute} < 10
\arabic{hour}:0\arabic{minute} on \number\month/\number\day/\number\year
\else
\arabic{hour}:\arabic{minute} on \number\month/\number\day/\number\year
\fi
}

\title{Convergence of invariant densities in the small-noise limit}

\author{Kevin K. Lin\\
{\bf klin@cims.nyu.edu}\\
Courant Institute of Mathematical Sciences, New York University\\
251 Mercer Street, New York, NY 10012}

\date{July 13, 2004}

\begin{document}

\maketitle

\begin{abstract}
\noindent
Let $\rho_0$ be an invariant probability density of a deterministic dynamical
system $f$ and $\rho_\eps$ the invariant probability density of a random
perturbation of $f$ by additive noise of amplitude $\eps.$ Suppose $\rho_0$ is
stochastically stable in the sense that $\rho_\eps\to\rho_0$ as $\eps\to 0.$
Through a systematic numerical study of concrete examples, I show that:
\begin{enumerate}

\item
The rate of convergence of $\rho_\eps$ to $\rho_0$ as $\eps\to 0$ is frequently
governed by power laws: $\onenorm{\rho_\eps-\rho_0}\sim\eps^\gamma$ for some
$\gamma>0.$

\item
When the deterministic system $f$ exhibits exponential decay of correlations, a
simple heuristic can correctly predict the exponent $\gamma$ based on the
structure of $\rho_0.$

\item
The heuristic fails for systems with some ``intermittency,'' {\em i.e.} systems
which do not exhibit exponential decay of correlations.  For these examples,
the convergence of $\rho_\eps$ to $\rho_0$ as $\eps\to 0$ continues to be
governed by power laws but the heuristic provides only an upper bound on the
power law exponent $\gamma.$

\end{enumerate}
Furthermore, this numerical study requires the computation of
$\onenorm{\rho_\eps-\rho_0}$ for $1.5 - 2.5$ decades of $\eps$ and provides an
opportunity to discuss and compare standard numerical methods for computing
invariant probability densities in some depth.
\end{abstract}


\section{Introduction}

Deterministic chaotic systems generally possess large numbers of invariant
probability measures.  Physical considerations lead naturally to the study of
invariant measures which are stable under small random perturbations.  (See
{\S{\ref{sec:bkgnd}}}.)  This paper explores the degree of stability of such
stochastically stable invariant measures in the setting of discrete-time
systems with invariant probability densities.

More precisely, consider dynamical systems of the form
\begin{equation}
\label{eqn:ds}
x_{k+1}=f(x_k)
\end{equation}
defined by a map $f:M\circlearrowleft,$ for example the map $x\mapsto 2x +
a\sin(2\pi x)\mbox{ (mod $1$)}$ shown in Figure {\ref{fig:map1a}}.  In this
paper, the space $M$ wil always be the circle $S^1,$ an interval in $\R,$ or a
product of such sets.  Thus operations like addition make sense in $M$ and
there is always a natural reference measure on $M,$ namely the normalized
Lebesgue measure.  Random perturbations of (\ref{eqn:ds}) can therefore be
defined via Markov chains of the form
\begin{equation}
\label{eqn:random-ds}
x_{k+1}=f(x_k)+\eps\xi_k,
\end{equation}
where the $(\xi_k)$ are independent, identically-distributed random variables
with a common probability density independent of $x_k.$ If $f$ is expanding
({\em i.e.} if the singular values of $Df(x)$ are all $>1$) for a sufficiently
large set of $x\in M$ and $M$ is connected, then the map $f$ possesses a unique
invariant probability density $\rho_0$ (see {\cite{baladi,baladi1}} for precise
statements of known results).  An invariant density $\rho_0$ is stochastically
stable if, for every continuous observable $\phi,$
$$\lim_{\eps\to 0}\int_M{\phi\rho_\eps} =\int_M{\phi\rho_0}$$ where $\int_M$
denotes integration with respect to Lebesgue measure on $M.$ When $\rho_0$ is
stochastically stable, it is natural to try to predict the rate with which
$\rho_\eps$ converges to $\rho_0$ as $\eps\to 0.$ As the examples in this paper
demonstrate, the rate of convergence of $\rho_\eps\to\rho_0$ can depend on the
detailed properties of $\rho_0$ and $f$ and may not always be immediately
apparent.

As a first step to understanding the factors which can affect the rate of
convergence in the limit $\eps\to 0,$ this paper systematically examines five
concrete examples numerically.  It is found that:
\begin{enumerate}

\item
The rate of convergence of $\rho_\eps$ to $\rho_0$ as $\eps\to 0$ is often
governed by power laws: $\onenorm{\rho_\eps-\rho_0}\sim\eps^\gamma$ for some
$\gamma>0,$ where $\onenorm{\rho_\eps-\rho_0} =\int_M{\abs{\rho_\eps-\rho_0}}$
denotes the $L^1$ norm on $M.$ Note that because
\begin{equation}
\label{eqn:onenorm}
\onenorm{\rho_\eps-\rho_0} =\sup_{\set{\phi:\norm{\phi}_\infty\leq 1}}
{\int_M{\phi\cdot(\rho_\eps-\rho_0)}},
\end{equation}
convergence in $L^1$ norm implies a uniform rate of convergence of expectation
values and is more stringent than stochastic stability.

\item
A simple heuristic allows one to predict the exponent $\gamma$ based on the
structure of the density $\rho_0,$ in a sense to be made more precise below
(see Equations (\ref{eqn:heuristic}) and (\ref{eqn:heuristic1})).

\item
The heuristic fails for systems with some ``intermittency,'' {\em i.e.} systems
which do not exhibit exponential decay of correlations
{\cite{liverani,pomeau}}.  For these examples, the convergence of $\rho_\eps$
to $\rho_0$ continues to be governed by power laws but the heuristic provides
only an upper bound on the power law exponent $\gamma.$

\end{enumerate}
The notation ``$f\sim g$'' means that there exist positive constants $c_1$ and
$c_2$ such that $c_1f(x)\leq g(x)\leq c_2f(x).$ Similarly, the notation
``$f\lesssim g$'' will be used below to mean that there exists $c>0$ such that
$f(x)\leq cg(x).$

Explaining the heuristic estimate requires the formalism of Perron-Frobenius
transfer operators {\cite{baladi}}.  Let $q$ be a function on $M$ and define
the operator $T_f$ by
\begin{equation}
\label{eqn:transop}
(T_fq)(x) = \sum_{y\in f^{-1}(x)} {\frac{q(y)}{\abs{\det(Df(y))}}}.
\end{equation}
The operator $T_f$ reformulates the dynamics in terms of probability densities:
if $q$ is the probability density of $x_k$, then $T_fq$ is the probability
density of $x_{k+1}=f(x_k).$ This elementary fact follows from the change of
variables formula.  Clearly, $\rho_0$ is an invariant density of $f$ if and
only if it is an eigenfunction of $T_f$ with eigenvalue 1: $T_f\rho_0=\rho_0.$
Furthermore, $T_f^nq\to\rho_0$ as $n\to\infty$ if the eigenvalue $1$ of $T_f$
is simple and the initial density $q$ is sufficiently regular.  And, if $T_f$
has a spectral gap, that is if the eigenvalue $1$ is an isolated point of the
spectrum $\sigma(T_f)$ of $T_f$ and the radius
$\abs{\sigma(T_f)\setminus\set{1}}$ of the smallest disc in $\C$ containing
$\sigma(T_f)\setminus\set{1}$ is strictly less than $1,$ then the convergence
of $T_f^nq$ to $\rho_0$ as $n\to\infty$ is exponentially fast.  The spectral
gap condition also implies the exponential decay of correlations, {\em i.e.}
there exist positive constants $c$ and $\theta<1$ such that
\begin{equation}
\abs{\int_M{\phi\cdot\of{T_f^n\psi}\cdot\rho_0}
  -\int_M{\phi\rho_0}\cdot\int_M{\psi\rho_0}}\leq c\theta^n
\end{equation}
for all sufficiently smooth observables $\phi$ and $\psi.$ Systems with
exponential decay of correlations are also said to be exponentially mixing.
Note that $1>\theta>\abs{\sigma(T_f)\setminus\set{1}}$.  See
{\cite{baladi,baladi1}} for details.

The ``noisy'' transfer operator is defined in a similar way: let $G_\eps$ be
the averaging operator
\begin{equation}
\label{eqn:avgop}
(G_\eps q)(x) =\eps^{-d}\int_M{g\of{\frac{x-y}{\eps}}q(y)\ dy},
\end{equation}
where $g$ is the common probability density of the IID random variables $\xi_k$
and $d=\dim(M).$ The operator $G_\eps$ represents the effect of additive noise
on a probability density $q.$ The noisy Perron-Frobenius operator $T_\eps$ is
then
\begin{displaymath}
T_\eps=G_\eps T_f.
\end{displaymath}
As before, the operator $T_\eps$ describes the random dynamics
(\ref{eqn:random-ds}) in terms of probability densities.  A density $\rho_\eps$
is invariant under (\ref{eqn:random-ds}) if and only if it is an eigenfunction
of $T_\eps$ of eigenvalue 1.

We can now discuss the heuristic.  Under fairly general conditions, $T_\eps^n
q$ converges to $\rho_\eps$ exponentially fast as $n\to\infty$ for fixed $\eps$
{\cite{baladi1}}.  Suppose $\rho_\eps\to\rho_0$ as $\eps\to 0.$ Since
$\lim_{n\to\infty}T_\eps^n\rho_0=\rho_\eps$ exponentially fast, we may expect
$\onenorm{T_\eps^n\rho_0-\rho_\eps}$ to be small for $n$ finite but
sufficiently large.  If such an $n$ can be chosen independent of $\eps,$ it is
then natural to guess that for all sufficiently small $\eps,$
\begin{equation}
\label{eqn:heuristic}
\onenorm{\rho_\eps-\rho_0}\sim\onenorm{T^n_\eps\rho_0-\rho_0}.
\end{equation}
Equation (\ref{eqn:heuristic}) states that one can estimate the difference
between $\rho_\eps$ and $\rho_0$ by studying the effect of applying $T_\eps$ a
finite number of times to the noiseless invariant density $\rho_0.$ Note that
the rate of convergence of $T_\eps^n\rho_0$ to $\rho_\eps$ as $n\to\infty$ is
governed by the size of the spectral gap
$\abs{\sigma(T_\eps)\setminus\set{1}}.$ If
$\abs{\sigma(T_\eps)\setminus\set{1}}$ can become arbitrarily small as $\eps\to
0,$ $n$ may have to be very large (or even increase as $\eps$ decreases) in
order for (\ref{eqn:heuristic}) to hold.  On the other hand, in exponentially
mixing systems with a sufficiently large spectral gap
$\abs{\sigma(T_\eps)\setminus\set{1}}$ ($\eps\geq 0$), it is likely that one
can take $n$ to be a small integer independent of $\eps.$ Equation
(\ref{eqn:heuristic}) should be taken only as a rough guideline for what can be
expected in studying the convergence of $\rho_\eps$ to $\rho_0$ in the $L^1$
norm in the small noise limit $\eps\to 0.$

In the exponentially mixing examples of Section {\ref{sec:gap}}, it is found
empirically that $n=1$ suffices:
\begin{equation}
\label{eqn:heuristic1}
\onenorm{\rho_\eps-\rho_0}\sim\onenorm{T_\eps\rho_0-\rho_0}
=\onenorm{G_\eps\rho_0-\rho_0}.
\end{equation}
This is not unexpected if
$\sup_{0\leq\eps<\eps_0} {\abs{\sigma(T_\eps)\setminus\set{1}}}\ll 1$ for some
$\eps_0>0.$ Equation (\ref{eqn:heuristic1}) states that one can estimate the
difference between $\rho_\eps$ and $\rho_0$ by studying the effect of applying
the averaging operator $G_\eps$ once to the noiseless invariant density
$\rho_0.$ Note that ``half'' of (\ref{eqn:heuristic1}) is always true: because
$(I-T_\eps)(\rho_\eps-\rho_0) = G_\eps\rho_0-\rho_0,$ the inequality
\begin{equation}
\label{eqn:half-heuristic}
\onenorm{G_\eps\rho_0-\rho_0}\leq 2\onenorm{\rho_\eps-\rho_0}
\end{equation}
always holds.  So as $\eps\to 0,$ $\onenorm{\rho_\eps-\rho_0}$ cannot converge
to $0$ faster than $\onenorm{G_\eps\rho_0-\rho_0}$ and the heuristic always
provides an upper bound on the power law exponent $\gamma.$

In problems where (\ref{eqn:heuristic1}) applies, the effect of $G_\eps$ on
$\rho_0$ (and hence the scaling of $\onenorm{G_\eps\rho_0-\rho_0}$ as $\eps\to
0$) depends on the structure of $\rho_0$: if $\rho_0$ is very smooth, as in the
uniformly expanding map in {\S\ref{sec:map1}}, then the effect of $G_\eps$ will
be to simply ``flatten'' $\rho_0,$ {\em i.e.} decrease the size of $\rho_0'$
(see Figure {\ref{fig:map1b}}).  The convergence of $G_\eps\rho_0$ to $\rho_0$
as $\eps\to 0$ will then be fast ({\em e.g.} $O(\eps^2)$), so the heuristic
also predicts a fast convergence of $\rho_\eps$ to $\rho_0$ as $\eps\to 0.$ On
the other hand, if $\rho_0$ contains any discontinuities or singularities (as
in the examples of {\S\ref{sec:map2}} -- {\S\ref{sec:map4}}), then the main
effect of $G_\eps$ will be to smooth out the discontinuities or singularities.
The rate of convergence predicted by the heuristic thus depends on the precise
form of the singularity in question.  As the heuristic is motivated by the
exponential convergence of $T_\eps^n\rho_0$ to $\rho_\eps$ (in the limit
$n\to\infty$) it is natural expect it to make sense only when
$\abs{\sigma(T_\eps)\setminus\set{1}} < 1$ uniformly in $\eps\geq 0.$ This
suggests, among other things, that the heuristic will work only when $T_f$
itself has a spectral gap, or that correlations decay exponentially fast in the
noiseless system (\ref{eqn:ds}).  As will be seen, the available numerical
evidence supports this claim.

\begin{table}
\begin{center}
\begin{tabular}{|c|c|c|}
\hline
& \multicolumn{2}{c|}{Exponent $\gamma=\lim_{\eps\to
      0}\frac{\log(\onenorm{\rho_\eps-\rho_0})}{\log(\eps)}$}\\\cline{2-3}
\raisebox{1.5ex}[0pt][0pt]{Map} & Predicted & Computed\\\hline\hline
Uniformly expanding & $2$ & {\mapone}\\\hline
Piecewise expanding & $1$ & {\maptwo}\\\hline
$\raisebox{-1.5ex}[0pt][0pt]{Quadratic (Misiurewicz)}$ &
\raisebox{-1.5ex}[0pt][0pt]{$\frac{1}{2}$} & {\mapthreea}\\
& & {\mapthreeb}\\\hline
& $\leq 0.5$ & {\mapfoura}\\
Neutral fixed point & $\leq 0.7$ & {\mapfourb}\\
& $\leq 0.3$ & {\mapfourc}\\\hline
Stadium & $\leq 2$ & {\mapfive}\\\hline
\end{tabular}
\caption{Summary of results.  The first example, a smooth uniformly expanding
  map, has a smooth invariant density.  Thus
  $G_\eps\rho_0-\rho_0=\frac{1}{6}\eps^2\rho_0''+o\of{\eps^2}$ and the
  heuristic predicts $O(\eps^2)$ convergence as $\eps\to 0.$ The second example
  has a piecewise continuous invariant density.  The heuristic thus predicts
  $O(\eps)$ convergence.  The third example has an invariant density which
  contain singularities of the form $x^{-1/2}$.  The heuristic thus predicts
  $\eps^{1/2}$ convergence.  The last two examples are not exponentially
  mixing, and the heuristic only predicts upper bounds.}
\label{tab:summary}
\end{center}
\end{table}

This paper examines five concrete examples, three of which are exponentially
mixing (see {\S\ref{sec:gap}}) and two which are not ({\S\ref{sec:nogap}}).  It
is found that the heuristic is valid for each of the exponentially mixing
examples but does not work for the intermittent examples.  Each example is
accompanied by an explanation of the numerical techniques used to compute the
invariant densities.  Table {\ref{tab:summary}} summarizes the numerical
results described in the rest of this paper.

The remainder of this section is devoted to a brief review of background and
motivation for this work.

\subsection{Motivation \& background}
\label{sec:bkgnd}

Chaotic systems are intrinsically unpredictable.  The tools of ergodic theory
and statistical mechanics are therefore necessary for questions concerning the
long-time behavior of dynamical systems {\cite{eckmann}}.  Invariant
probability measures capture the long-time statistical properties of dynamical
systems.

Deterministic dynamical systems generally possess a multitude of invariant
probability measures.  For example, fixed points and periodic orbits support
invariant measures.  However, when the system in question is chaotic, most of
these invariant measures are unstable and do not have directly observable
effects on the long-time dynamics.  The notion of stochastic stability (due to
Kolmogorov) arises naturally as a criterion for separating these unstable
invariant measures from those which are physically relevant: because no real
physical system can be completely isolated from the rest of the universe, every
physical experiment is susceptible to the effects of noise.  Thus only
stochastically stable invariant measures have directly observable effects on
the long-term statistics of dynamical systems.

One of the earliest general results on the stochastic stability of invariant
measures concerns systems with uniformly hyperbolic attractors.  Such systems
possess Sinai-Ruelle-Bowen (SRB) measures {\cite{young}}, which (among their
many properties) capture the statistical properties of a set of initial
conditions of positive Lebesgue measure.  Kifer proved that SRB measures of
uniformly hyperbolic systems are stochastically stable {\cite{kifer}}.
However, his method of proof cannot address the question of convergence rate in
the limit of small noise.  In {\cite{baladi1}}, Baladi and Young examine the
question of stochastic stability and the rate of convergence in the setting of
expanding maps and convolution-type perturbations.  Blank and Keller, in
subsequent work, proved results for more general perturbations of
$1$-dimensional maps {\cite{blank1}}.

In addition to understanding more deeply the stochastic stability of invariant
measures, the results described here may be relevant for numerical studies of
dynamical systems with intermittent, metastable behavior.  The injection of
noise into a numerical simulation can help reduce initialization bias in
numerical computations of expectation values via time averaging while adding
controllable errors to the computed expectation value {\cite{lin1}}.


\section{Exponentially mixing systems}
\label{sec:gap}

The examples in this section have one common feature: their transfer operators
$T_f$ all have spectral gaps.  That is, $\abs{\sigma(T_f)\setminus\set{1}}<1.$
This implies the exponential decay of correlations.  The heuristic estimate
(\ref{eqn:heuristic}) or its variant (\ref{eqn:heuristic1}) are expected to
work for maps whose transfer operators have spectral gaps.  This is found to be
the case.  For all the examples in this section, the common density $g:\R\to\R$
of the random variables $\xi_k$ is taken to be
\begin{equation}
\label{eqn:pert}
g(x)=\left\{\begin{array}{ll}
\frac{1}{2},&\abs{x}\leq 1\\
0,&\abs{x}> 1\\
\end{array}\right..
\end{equation}
(The first two examples are maps on $S^1,$ here identified with the interval
$[0,1].$ The random perturbations should therefore be taken modulo $1.$) The
perturbations $\eps\xi_k$ are thus uniform random variables on the interval
$[-\eps,+\eps].$ Calculations using gaussian kernels wrapped around the circle
(not shown here) indicate that the scaling is insensitive to the exact form of
the kernel $g.$

\subsection{Smooth expanding circle map}
\label{sec:map1}

\begin{figure}
\begin{center}
\putgraph{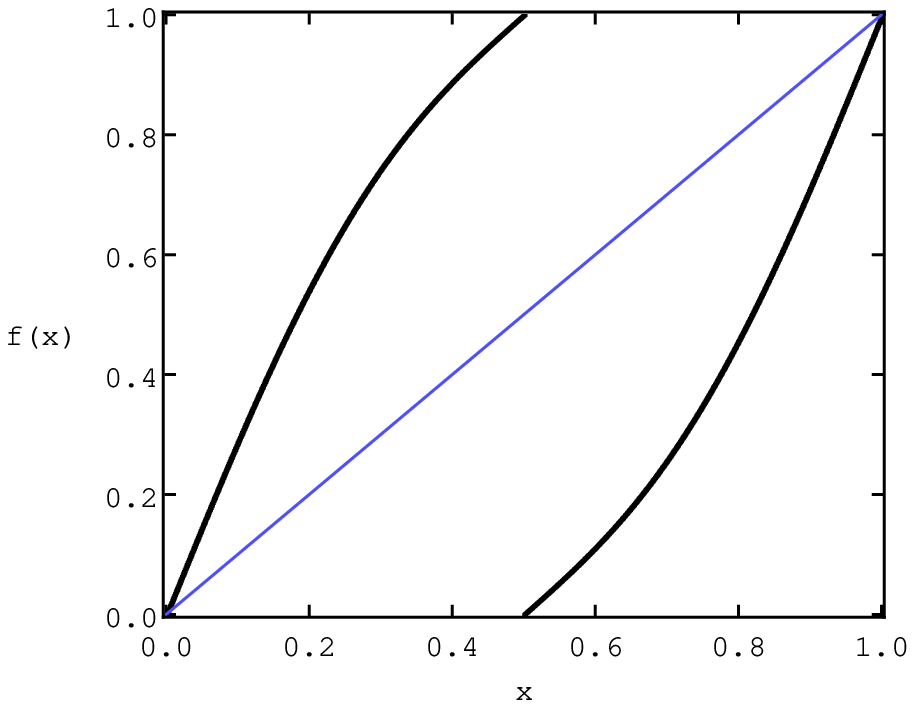}
\caption{A smooth uniformly expanding map.}
\label{fig:map1a}
\end{center}
\end{figure}

\begin{figure}
\begin{center}
\putgraph{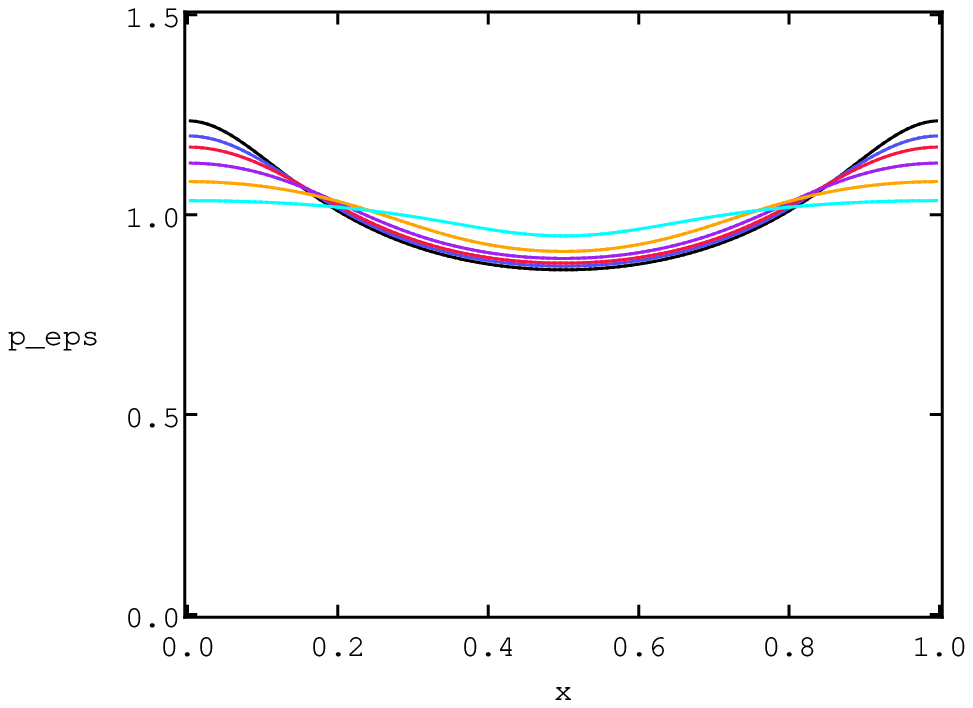}
\caption{The invariant densities $\rho_\eps$ with $\eps\in\set{0, \frac{1}{10},
    \frac{1}{10\sqrt{2}}, \frac{1}{5}, \frac{1}{5\sqrt{2}}, \frac{2}{5}},$ for
    the map (\ref{eqn:map1}).  The black curve is $\rho_0$; the densities
    $\rho_\eps$ become flatter as $\eps$ increases.  The map parameter is
    $a=0.15.$}
\label{fig:map1b}
\end{center}
\end{figure}

Consider $f:S^1\circlearrowleft$ defined by
\begin{equation}
\label{eqn:map1}
f(x)=2x + a\sin(2\pi x)\mbox{ (mod 1)}.
\end{equation}
See Figure {\ref{fig:map1a}}.  This map is clearly smooth and is uniformly
expanding for $0 < a <\frac{1}{2\pi}.$ Its invariant density $\rho_0$ is easy
to compute by discretizing the Perron-Frobenius operator $T_f$ on a uniform
mesh and is shown in Figure {\ref{fig:map1b}}.

Before discussing the details of numerical calculations, it should be noted
that one can prove (\ref{eqn:heuristic1}) for the map (\ref{eqn:map1}) when the
parameter $a$ is sufficiently small.  Although the analysis is simple and uses
only standard techniques, it is nevertheless instructive and is included here
for completeness.  Let $B$ denote the space $C^1$ with the usual norm
$\norm{h}_B=\max\of{\norm{h}_\infty,\norm{h'}_\infty},$ and let $B_0$ be the
subspace $\set{h\in B:\int_M{h}=0}.$ Note that $B_0$ is $T_\eps$-invariant for
$\eps\geq 0.$ It is not difficult to show that the Perron-Frobenius operator
$T_f$ associated with (\ref{eqn:map1}) has the property that
$\norm{T_f}_{B_0}<1$ when the parameter $a$ is sufficiently small:
differentiating (\ref{eqn:transop}) shows that $\norm{(T_fq)'}_\infty\leq
c_a\norm{q}_B$ for some constant $c_a$ dependent on $a,$ and $c_a<1$ when $a$
is small.  Because $\int_M{T_fq}=\int_M{q}=0,$ the Poincar\'e inequality shows
that $\norm{T_fq}_\infty\leq\norm{(T_fq)'}_\infty.$ So $\norm{T_f}_{B_0}
=\sup\set{\norm{T_fq}_B:q\in B_0,\norm{q}_B=1}\leq c < 1.$

Now let $g$ be a probability density on $S^1$ and let $G_\eps$ be the
associated averaging operator.  Denote the (unique) invariant density of
$T_\eps =G_\eps T_f$ by $\rho_\eps.$ Since $\norm{G_\eps}_B\leq 1,$
$\norm{T_\eps}_{B_0}\leq\norm{T_f}_{B_0}<1$ for all $\eps\geq 0.$ Thus the
restriction $\restrict{(I-T_\eps)}{B_0}$ of $I-T_\eps$ to $B_0$ has bounded
inverse for all $\eps\geq 0.$ A simple calculation shows that
\begin{equation}
\rho_\eps-\rho_0
=\brac{\restrict{(I-T_\eps)}{B_0}}^{-1}\of{G_\eps\rho_0-\rho_0}.
\end{equation}
Thus
\begin{align*}
\norm{\rho_\eps-\rho_0}_B
&\leq
\norm{(I-T_\eps)^{-1}}_{B_0}\norm{G_\eps\rho_0-\rho_0}_B\\
&\leq\of{1-\norm{T_f}_{B_0}}^{-1}\norm{G_\eps\rho_0-\rho_0}_B.
\end{align*}
Also
\begin{align*}
\norm{G_\eps\rho_0-\rho_0}_B
&=\norm{(I-T_\eps)\cdot(I-T_\eps)^{-1}\cdot(G_\eps\rho_0-\rho_0)}_B\\
&\leq\norm{I-T_\eps}_{B_0}\cdot\norm{\rho_\eps-\rho_0}_B.
\end{align*}
Since $\norm{I-T_\eps}_{B_0}\leq 1 +\norm{T_f}_{B_0}$, this means
\begin{equation}
\norm{\rho_\eps-\rho_0}_B\sim\norm{G_\eps\rho_0-\rho_0}_B
\end{equation}
with $\eps$-independent constants in the $\sim$ relation.  Note that this
already means $\rho_\eps-\rho_0$ converges to $0$ at the same rate as
$G_\eps\rho_0-\rho_0$ in the $C^1$ metric as $\eps\to 0.$ Combining Equation
(\ref{eqn:half-heuristic}) and the fact that
$\norm{\rho_\eps-\rho_0}_B\geq\onenorm{\rho_\eps-\rho_0}$ yields
\begin{equation}
\onenorm{G_\eps\rho_0 -\rho_0}\lesssim\onenorm{\rho_\eps
  -\rho_0}\lesssim\norm{G_\eps\rho_0 -\rho_0}_B.
\end{equation}
For the map (\ref{eqn:map1}), it can be shown that the unique invariant density
$\rho_0$ is smooth (see {\cite{baladi}} for details).  So $G_\eps\rho_0-\rho_0
=\frac{1}{6}\eps^2\rho_0''+o\of{\eps^2}$ and
$\onenorm{G_\eps\rho_0-\rho_0}\sim\norm{G_\eps\rho_0-\rho_0}_B$ when $\eps$ is
sufficiently small, so that
\begin{equation}
\onenorm{G_\eps\rho_0-\rho_0}\sim\onenorm{\rho_\eps-\rho_0}\sim\eps^2.
\end{equation}
Thus the heurstic estimate holds for (\ref{eqn:map1}) when $a$ is small.  A
little more is true: for any smooth expanding map $f,$ $1$ is an eigenvalue of
the Perron-Frobenius operator $T_f$ with algebraic multiplicity $1$ (see
{\cite{baladi}}).  This means $\restrict{T_f}{B_0}$ has spectral radius $<1.$
As $\abs{\spec{T_f}}=\lim_{n\to\infty}\norm{T_f^n}^{1/n},$ the facts above
together with the previous argument show that for any smooth expanding map $f$
there exists an integer $N\geq 1$ such that the heuristic holds for $f^N=f\circ
f\circ ...\circ f$, {\em i.e.} if $\rho^{(N)}_\eps$ is the invariant density
for $x_{k+1}=f^N(x_k)+\eps\xi_k$ then
$\onenorm{\rho^{(N)}_\eps-\rho_0}\sim\onenorm{G_\eps\rho_0-\rho_0}.$

To test the general validity of the heuristic (\ref{eqn:heuristic1}) when $a$
is not necessarily small, it is natural to choose $a$ close to $1/2\pi =
0.15915494309189535....$ In what follows, $a$ is set to $0.15.$ In Figure
{\ref{fig:map1b}}, it can be seen that $\rho_\eps$ becomes flatter as $\eps$
increases.  This is not suprising: as $\eps\to\infty$, the noise dominates the
dynamics and the invariant density is just the Lebesgue measure.  As $\eps\to
0,$ the densities converge to $\rho_0.$

\begin{figure}
\begin{center}
\putgraph{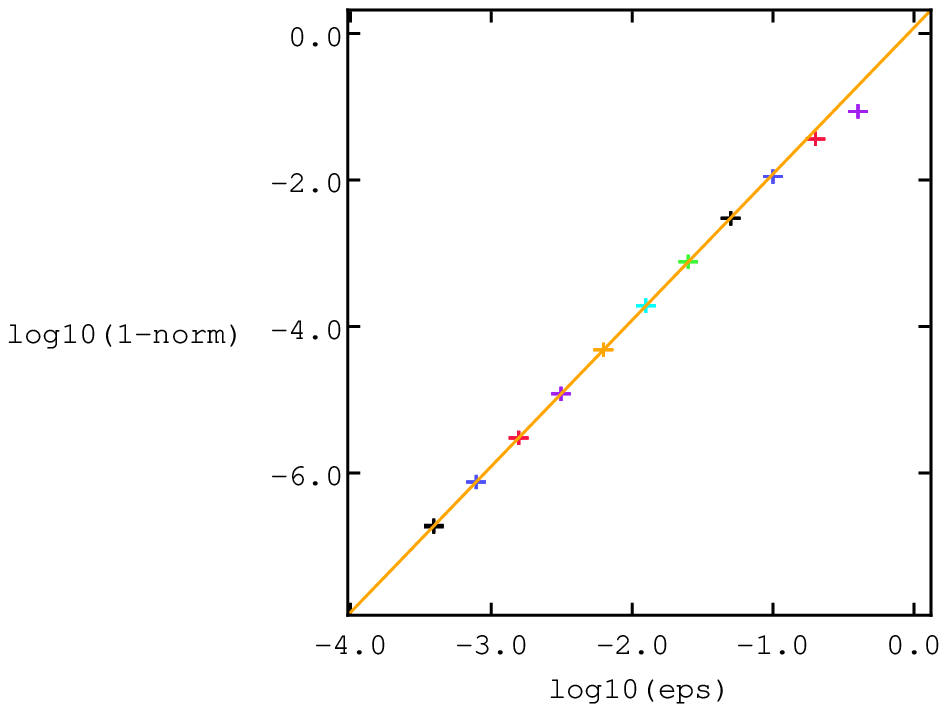}
\caption{The differences $\onenorm{\rho_\eps-\rho_0}$ and
  $\onenorm{G_\eps\rho_0-\rho_0}$ as a function of $\eps$ on a log-log graph.
  The slope is {\mapone}. The slope is calculated using least-squares
  regression on the interval $\eps\in\maponerange.$ Note that this means data
  points on the far right are discarded.  Throughout the paper, data points for
  which error bars are available are drawn as vertical lines with 3 horizontal
  marks: the vertical line marks the abscissa of the data point, the middle
  mark the ordinate, and the top and bottom marks are the upper and lower
  bounds on the error.  The error bars in this figure are too small to be
  seen.}
\label{fig:map1c}
\end{center}
\end{figure}

\begin{figure}
\begin{center}
\putgraph{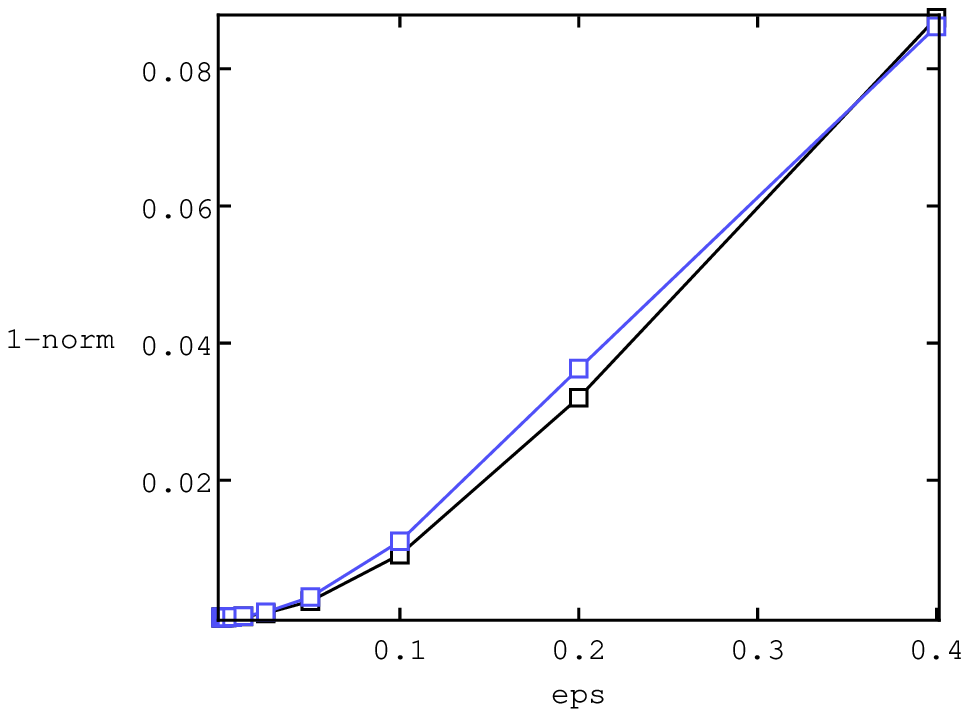}
\caption{The differences $\onenorm{\rho_\eps-\rho_0}$ and
  $\onenorm{G_\eps\rho_0-\rho_0}$ as a function of $\eps$ on a linear graph.}
\label{fig:map1d}
\end{center}
\end{figure}

The $L^1$ distances $\onenorm{\rho_\eps-\rho_0}$ are plotted in Figure
{\ref{fig:map1c}} as a function of $\eps$ on a log-log scale.  Clearly,
$\onenorm{\rho_\eps-\rho_0}\sim\eps^\gamma$ for some $\gamma.$ The exponent
$\gamma$ can be calculated as the slope of the line in Figure {\ref{fig:map1c}}
using standard least-squares regression: it is about {\mapone} when $\gamma$ is
computed via least-squares regression on the interval $\eps\in\maponerange.$

In comparison, $G_\eps\rho_0 -\rho_0 =\frac{1}{6}\eps^2\rho_0''+o\of{\eps^2}$
implies $\onenorm{G_\eps\rho_0 -\rho_0}\sim\eps^2.$ Thus, the heuristic
estimate (\ref{eqn:heuristic1}) predicts that
$$\onenorm{\rho_\eps-\rho_0}\sim\eps^2.$$ For this particular map, the
invariant density $\rho_0$ is smooth.  Therefore the main effect of $G_\eps$ on
$\rho_0$ is to ``flatten'' it out, consistent with the pictures of $\rho_\eps$
in Figure {\ref{fig:map1b}}.

The numerical data thus provide evidence that the heuristic is valid for values
of the parameter $a$ which is not close to $0.$ The heuristic argument predicts
a little bit more than just the exponent $\gamma$: it predicts that when the
spectral gap of $T_f$ is sufficiently large, the quantities
$\onenorm{\rho_\eps-\rho_0}$ and $\onenorm{G_\eps\rho-\rho_0}$ will be
comparable.  This is also borne out by the data, as can be seen in Figure
{\ref{fig:map1d}}.  The numerically computed spectral gap is about $0.848,$
which is quite large.  This explains why $\onenorm{\rho_\eps-\rho_0}$ and
$\onenorm{G_\eps\rho_0-\rho_0}$ are so close in Figure {\ref{fig:map1d}}.

\subsubsection*{Numerical method.}

It is straightforward to compute the invariant densities shown in Figure
{\ref{fig:map1b}} using a finite difference discretization of the
Perron-Frobenius operator $T_f.$ Let $\phi_j$ denote the inverse branches of
$f$: $f\circ\phi_j=\mbox{id}, j=1,2.$ (There are two branches because the map
$f$ is $2$-to-$1$.)  We can then write $T_f$ as $$(T_fq)(x)
=\sum_{j=1}^{2}{\frac{q(\phi_j(x))}{\abs{f'(\phi_j(x))}}}.$$ The numerical
procedure is then:
\begin{enumerate}

\item
Identify $S^1$ with the interval $[0,1]$ with periodic boundary conditions and
let $\hat{x}_i=i/N, i=0,1,2,...,N-1,$ be a uniform grid of size $N.$

\item
Let $\hat{\rho}$ be an $N$-vector with nonnegative entries.  Define the matrix
$\hat{T}$ by
\begin{equation}
(\hat{T}\hat{\rho})_i
  =\sum_{j=1}^{2}{\frac{\hat{\rho}(\phi_j(\hat{x}_i))}
  {\abs{f'(\phi_j(\hat{x}_i))}}},
\end{equation}
where we extend the $N$-vector $\hat{\rho}$ to a function on $[0,1]$ via
polynomial interpolation using the grid points closest to $x.$ Note that
interpolation is necessary becuase the vector $\hat\rho$ only contains values
of the density on the grid and $\phi_j(\hat{x}_i)$ will generally not be a grid
point.

\item
Using a numerical linear algebra package like ARPACK {\cite{arpack}}, which
implements an iterative Arnoldi eigenproblem solver, compute the eigenvector
$\hat\rho_0$ of $\hat{T}$ with the largest eigenvalue.  This requires only a
function for performing matrix-vector multiplies, which avoids storing and
diagonalizing the entire matrix $\hat{T}.$

\end{enumerate}
In the computations shown in this section, the polynomial interpolation of
$\hat\rho$ utilizes a stencil with $6$ points and the grid consists of $N=10^4$
points.  For smooth maps like (\ref{eqn:map1}), the order of convergence (as
$N\to\infty$) is formally $6.$ This is confirmed empirically by numerical
tests.  Thus the finite difference scheme converges rather rapidly and provides
an efficient way to compute $\rho_0.$

Note that there is a closely related method, due to Ulam {\cite{keane}}, which
consists of forming a partition of the phase space and computing transition
probabilities to form a stochastic matrix.  The stochastic matrix corresponds
to a finite-state Markov chain which mimics the dynamics of the map $f$ and
whose invariant probability distribution can be computed using numerical
eigenproblem solvers.  Ulam's method is more general: it applies to dynamical
systems whose invariant measures are singular.  As a method for computing
invariant densities, however, it is empirically found not to be as accurate or
efficient as the high-order finite difference schemes outlined above.  This is
perhaps not so surprising: Ulam's method is, intuitively, only a low-order
finite element scheme for the equation $T_f\rho_0 =\rho_0.$

The computation of $\rho_\eps$ relies on a similar finite difference scheme.
Instead of discretizing $T_\eps = G_\eps T_f$ directly, it is easier to reuse
the matrix $\hat{T}$ from the computation of $\rho_0$ and multiply it from the
left by a discretization $\hat{G}$ of $G_\eps.$ The matrix $\hat{G}$ is
constructed by
\begin{equation}
(\hat{G}\hat{\rho})_i
  =\frac{1}{2\eps}\int_{\hat{x}_i-\eps}^{\hat{x}_i+\eps}{\hat{\rho}(x)\ dx}.
\end{equation}
The integral is evaluated numerically using the trapezoid rule.  This simple
scheme suffices here because it is second-order accurate, and with $N=2\times
10^4$ the expected error is on the order of $10^{-8}.$ This is smaller than the
$L^1$ norms $\onenorm{\rho_\eps-\rho_0}$ to be computed.  ARPACK can then be
applied to the matrix $\hat{G}\cdot\hat{T}.$ The reason that $T_f$ is
discretized to $6$th order while $G_\eps$ is only discretized to $2$nd order is
empirical: numerical experiments indicate that discretizing $T_f$ to higher
order accelerates the convergence of ARPACK eigenproblem routines.  This
phenomenon is not completely understood, but intuitively a higher-order
discretization $\hat{T}$ of $T_f$ reproduces the spectral properties of $T_f$
more accurately.  As $\restrict{T_f}{B_0}$ is strongly contractive, $\hat{T}$
will also be strongly contractive.

In actual computations performed on Sun Ultra 5 workstations and on PowerPC
G3/G4-based Macintosh computers, the computation of each invariant density
$\rho_\eps$ can take anywhere from $10$ seconds to $5$ minutes.  The variations
in running time are sometimes quite unpredictable; they are certainly not
monotic in $\eps$ as one might expect.  Note that the running time of ARPACK
routines depends quite sensitively on the spectral decomposition of $\hat{T},$
so the variations in running times provide some indirect evidence in the
complexity of the $\eps$-dependence of the spectrum of $T_\eps.$ In any case,
the finite difference scheme is in general quite efficient and allows the
systematic exploration of different values of $\eps$ (shown above) and $a$ (not
shown here).

The error bars in Figure {\ref{fig:map1c}} are computed by repeating the
calculations with $N=10^4$ points and checking the numerical convergence of the
resulting estimates for $\onenorm{\rho_\eps-\rho_0}.$ This provides error
estimates $\Delta A(\eps)$ for each computed value of
$A(\eps)=\log_{10}(\onenorm{\rho_\eps-\rho_0}).$ These error estimates can then
be combined to provide an upper bound on the error $\Delta\gamma$ in the
estimated exponent $\gamma$ by linearizing $\gamma$ about the computed values
$A(\eps_i)$:
\begin{equation}
\Delta\gamma
\lesssim\frac{1}{n}\sum_{i=1}^{n}{\abs{\frac{\log_{10}(\eps_i)-\mu}
{\sigma^2}\cdot\Delta A(\eps_i)}},
\end{equation}
where $\mu=\frac{1}{n}\sum_{i=1}^{n}{\log_{10}(\eps_i)}$ and
$\sigma^2=\frac{1}{n}\sum_{i=1}^{n}{(\log_{10}(\eps_i)-\mu)^2}.$ Note that
because of the high order of convergence of the finite difference scheme, the
main source of error in estimating $\onenorm{rho_\eps-\rho_0}$ actually comes
from the evaluation of $\onenorm{\hat\rho_\eps-\hat\rho_0}$ using the trapezoid
rule.  Since $\rho_\eps-\rho_0$ is continuous, the trapezoid rule is
second-order, and this error dominates all others in the computation.  It is
therefore natural to combine the results of the two runs to obtain a more
accurate answer via Richardson extrapolation.  This has been done in computing
the exponent $\gamma.$ The error estimates are therefore rather conservative.


\subsection{Piecewise-expanding map}
\label{sec:map2}

The second example is
\begin{equation}
\label{eqn:map2}
f(x)=\left\{\begin{array}{ll}
x+\frac{1}{2},& 0\leq x\leq\frac{1}{2}\\
2(1-x),       & \frac{1}{2}< x\leq 1\\
\end{array}\right..
\end{equation}
Strictly speaking, $f$ is not expanding because $f'(x)=1$ for $0\leq
x\leq\frac{1}{2}.$ However, $f^2=f\circ f$ is expanding and $f$ falls
within the class of piecewise expanding maps considered by Lasota and
Yorke {\cite{lasota}}.

\begin{figure}
\begin{center}
\putgraph{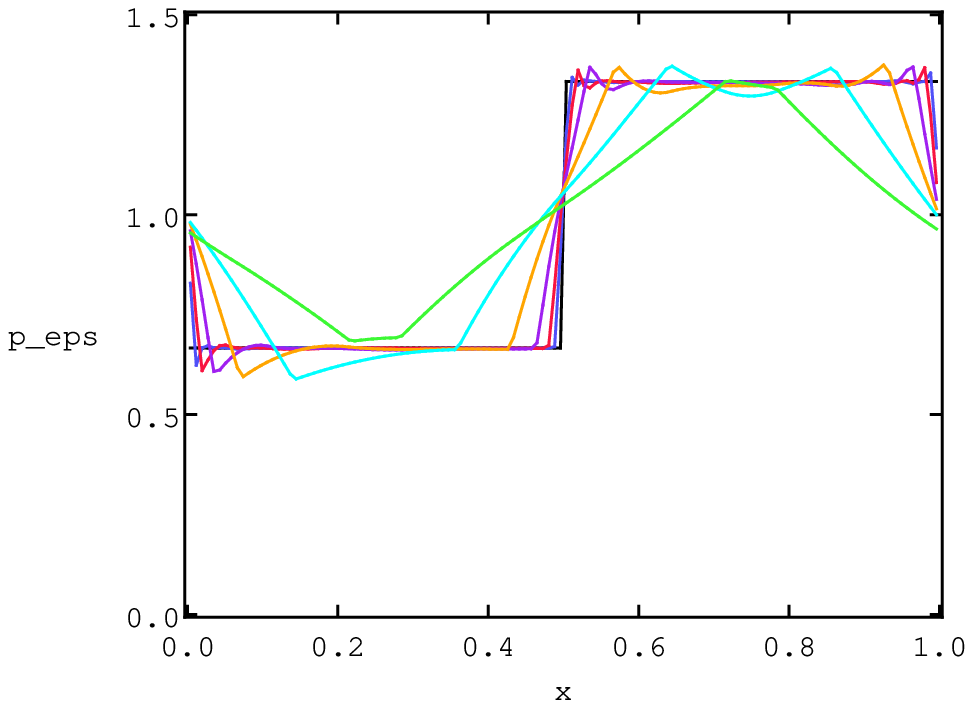}
\caption{Invariant densities $\rho_\eps$ for (\ref{eqn:map2}) and its random
  perturbations, with $\eps$ = 0 (black curve), 0.0088, 0.0176, 0.0353, 0.0707,
  0.141, and 0.282.  As $\eps$ increases, the width of the transition region
  increases too.}
\label{fig:map2b}
\end{center}
\end{figure}

The invariant density of (\ref{eqn:map2}) is easy to compute analytically: it
is
\begin{displaymath}
\rho_0(x)=\left\{\begin{array}{ll}
\frac{2}{3},& 0\leq x\leq\frac{1}{2}\\
\frac{4}{3},& \frac{1}{2}< x\leq 1\\
\end{array}\right..
\end{displaymath}
The invariant densities $\rho_\eps$ are plotted in Figure {\ref{fig:map2b}}.
Again, it can be seen that $\rho_\eps$ becomes flatter as $\eps$ increases: the
main effect of $G_\eps$ is to smear out the jump discontinuities of $\rho_0.$
Note that the ``overshoot'' visible in Figure {\ref{fig:map2b}} at $x=0$ and
$x=\frac{1}{2},$ reminiscent of the Gibbs phenomenon, is not a numerical
artifact.  Rather it is the product of the action of the map and the random
perturbation.  More precisely, it is produced by the following mechanism:
\begin{enumerate}

\item
Consider $T_\eps^n\rho_0$ with $n=1,2,3,...$.  For $n=1,$ $T_\eps\rho_0 =
G_\eps\rho_0$ looks very much like $\rho_0,$ but the action of $G_\eps$
replaces the discontinuity with a linear transition region of $O(\eps)$ width.

\item
On the next iterate, the action of $T_f$ then cuts $G_\eps\rho_0$ into two
pieces and rearranges them in such a way as to produce the ``overshoot'' of
$O(1)$ height and $O(\eps)$ width.  The next application of $G_\eps$ smooths
out the transition some more but does not change the asymptotic
$\eps$-dependence of the overshoot and the transition region.

\end{enumerate}
As $T_f^n\rho_0$ converges to $\rho_\eps$ with increasing $n,$ this process is
iterated to produce the overshoot structure.  Seen in this light, the overshoot
must have $O(1)$ height and $O(\eps)$ width, making a $O(\eps)$ contribution to
$\onenorm{\rho_\eps-\rho_0}.$

\begin{figure}
\begin{center}
\putgraph{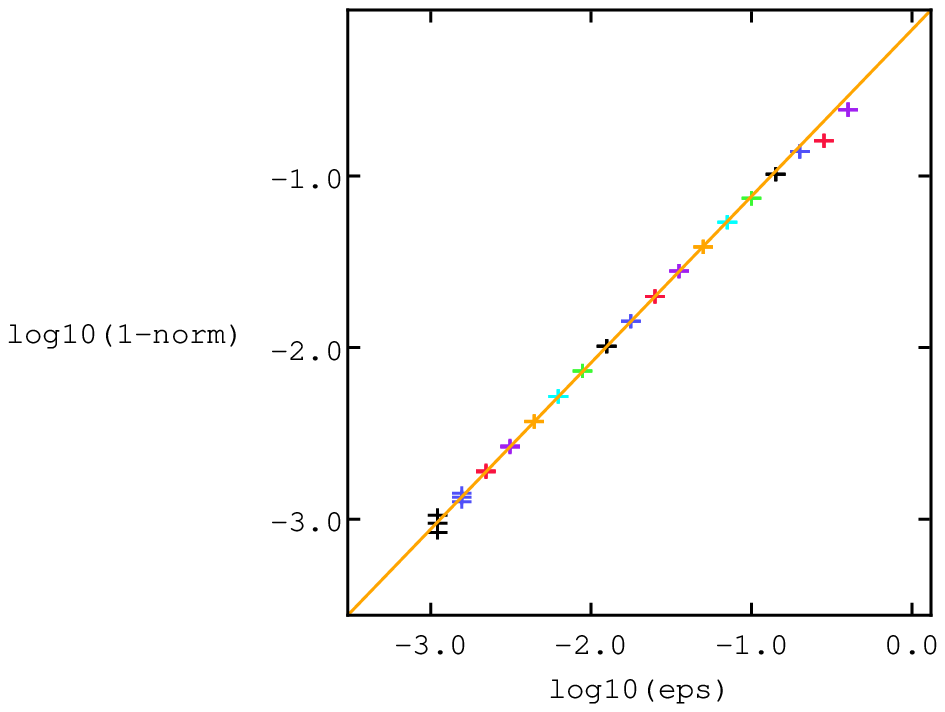}
\caption{The differences $\onenorm{\rho_\eps-\rho_0}$ and
  $\onenorm{G_\eps\rho_0-\rho_0}$ as a function of $\eps$ on a log-log graph.
  The slope is {\maptwo}.  The slope is calculated using least-squares
  regression on the interval $\eps\in\maptworange.$ Note that some of the data
  points have error bars which are too small to be seen on this scale.}
\label{fig:map2c}
\end{center}
\end{figure}

\begin{figure}
\begin{center}
\putgraph{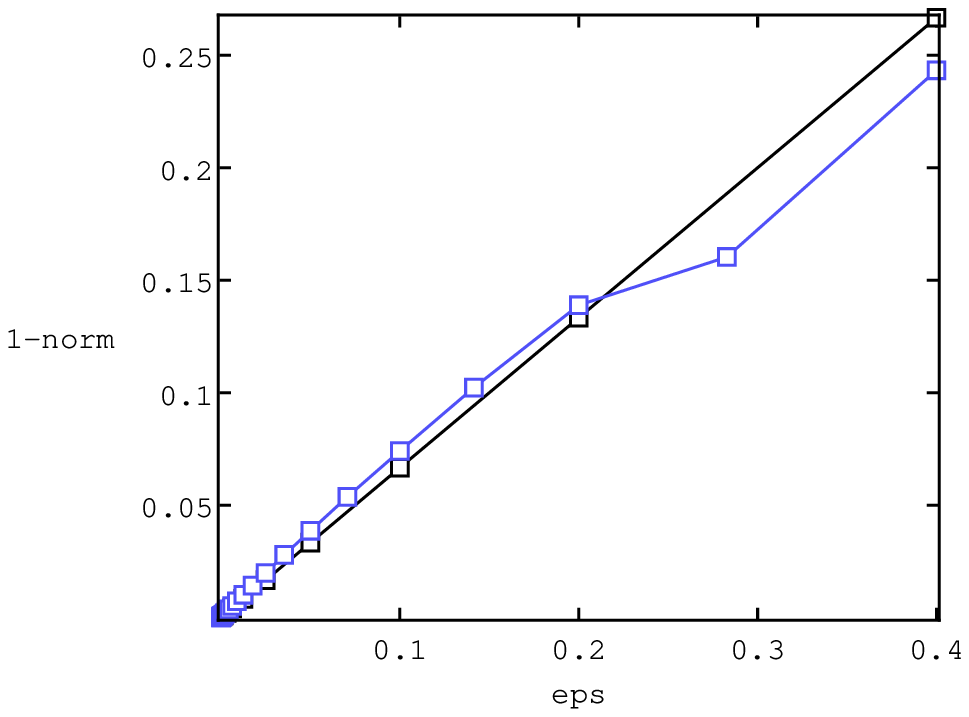}
\caption{The differences $\onenorm{\rho_\eps-\rho_0}$ and
  $\onenorm{G_\eps\rho_0-\rho_0}$ as a function of $\eps$ on a linear graph.}
\label{fig:map2d}
\end{center}
\end{figure}

Clearly, the main effect of $G_\eps$ is to smear out the discontinuity in
$\rho_0.$ So $\onenorm{G_\eps\rho_0-\rho_0}\sim\eps$ and the heuristic predicts
\begin{displaymath}
\onenorm{\rho_\eps-\rho_0}\sim\eps.
\end{displaymath}
This is consistent with the densities $\rho_\eps$ shown in Figure
{\ref{fig:map2b}}.  These $L^1$ norms are plotted in Figure {\ref{fig:map2c}}
as a function of $\eps$ on a log-log scale: the slope $\gamma$ is {\maptwo},
consistent with the prediction that $\gamma=1.$ The linear scaling of
$\onenorm{\rho_\eps-\rho_0}$ with $\eps$ should also be apparent on a linear
scale.  This is indeed the case: see Figure {\ref{fig:map2d}}.

The proof from Section {\ref{sec:map1}} can be adapted to show that there
exists an $N>0$ such that if $\rho^{(N)}_\eps$ is the unique invariant density
of $x_{k+1} = f^N(x_k)+\eps\xi_k,$ then $T_{f^N}$ is a bounded operator on the
space $B$ of functions of bounded total variation with the norm
$\norm{\cdot}_{BV} =\onenorm{\cdot}+\mbox{(total variation)}$ and
$\norm{T_{f^N}}_{B_0}<1.$ One obtains
$\norm{\rho_\eps-\rho_0}_B\sim\norm{G_\eps\rho_0-\rho_0}_B.$ This result is,
unfortunately, not so useful for this map: because $\rho_0$ is discontinuous,
$G_\eps\rho_0$ will not converge to $\rho_0$ in $\norm{\cdot}_B$ as $\eps\to
0.$

\subsubsection*{Numerical method.}

For this example, it is straightforward to implement the finite difference
scheme of Section {\ref{sec:map1}}, with the following modifications:
\begin{enumerate}

\item
One should avoid interpolating across discontinuities in $f$ and in $f'.$ Doing
so is empirically seen to produce unacceptably large errors in the computation
of $\hat{T}\hat{\rho}.$ For (\ref{eqn:map2}), this means the polynomial
interpolation of $\hat{\rho}$ should not use stencils which contain $0$ or
$\frac{1}{2}$ as an interior point.

\item
In order to avoid interpolating across $0$ or $\frac{1}{2},$ it is necessary to
do one-sided interpolations, {\em i.e.} interpolate near the edge of the
stencil.  It is well known that one-sided polynomial interpolation on uniform
grids can produce large errors.  Rather than using Legendre interpolation,
however, it suffices to simply increase the number of grid points and check
that the resulting answer has converged numerically.

\end{enumerate}
In the computations shown in this section, the polynomial interpolation of
$\hat{\rho}$ utilizes a stencil with $4$ points and the grid again consists of
$N=10^4$ points.  Numerical tests show that the order of accuracy lies between
$3$ and $4.$ Thus the finite difference scheme is still sufficiently accurate
to provide an efficient way to compute $\rho_\eps.$

The error bars in Figure {\ref{fig:map2c}} are computed by repeating the
calculations with $N=5\times 10^3$ points and checking the numerical
convergence of the resulting estimates for $\onenorm{\rho_\eps-\rho_0}.$ As in
{\S\ref{sec:map1}}, the main source of error in estimating
$\onenorm{\rho_\eps-\rho_0}$ comes from the trapezoid rule.  The grid
explicitly contains the points of discontinuity, so the trapezoid rule is still
second-order accurate.


\subsection{Quadratic maps}
\label{sec:map3}

The third example is the quadratic map $f:[-1,1]\circlearrowleft$:
\begin{equation}
\label{eqn:map3}
f(x)=0.9-ax^2, a > 0.
\end{equation}
Note that $f$ maps $[-1,1]$ into $[0.9-a,0.9],$ so for $a\leq 1.8$ and
$\eps\leq 0.1,$ the Markov chain defined by (\ref{eqn:random-ds}) will always
stay inside the interval $[-1,1].$

Unlike the previous examples, this map has a critical point: $f'(0)=0.$ The map
is contractive in a neighborhood of this critical point.  Whether $f$ has a
invariant density (as opposed to singular invariant measures) depends on the
fate of the critical point.  For example, if there exists an integer $n>0$ and
a stable fixed point $x_0$ such that $f^n(0)=x_0,$ then a positive amount of
probability will collapse onto $x_0$ to create a $\delta$ mass there.  Another
scenario which can prevent the existence of an invariant density is if $f^n(0)$
comes close to $0$ infinitely often as $n\to\infty.$

\begin{figure}
\begin{center}
\putgraph{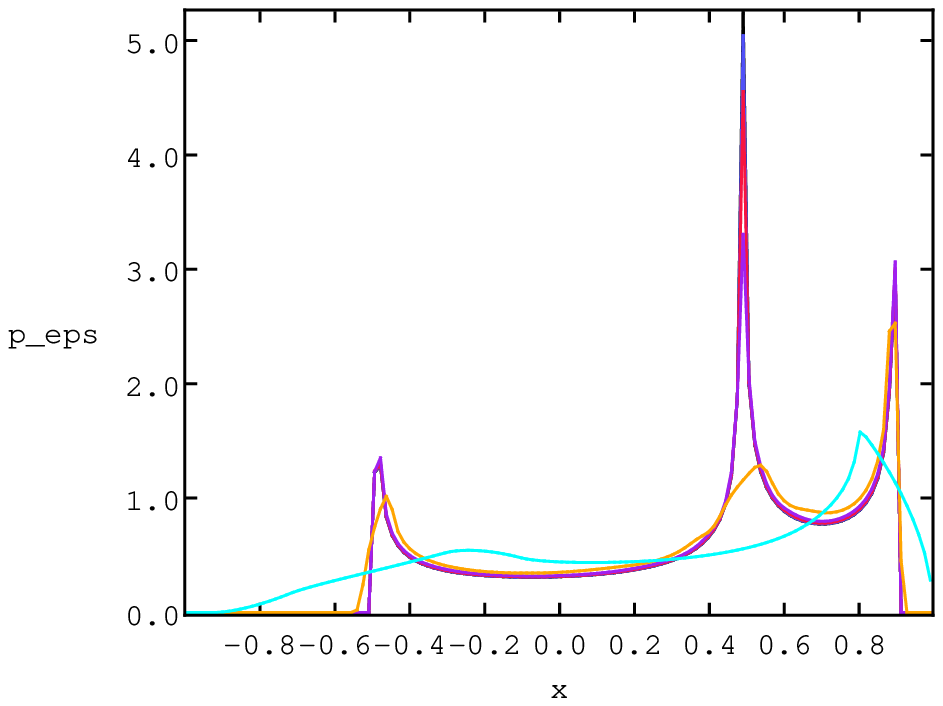}
\caption{The invariant densities $\rho_\eps$ with increasing values of
  $\eps.$ The map parameter is $a=1.7152100141023...$; the critical
  point is sent to an unstable fixed point at $0.489320111422868...$ in
  $3$ steps.  The noiseless density therefore has three $x^{-1/2}$
  singularities.}
\label{fig:map3b}
\end{center}
\end{figure}

To ensure the existence of a density, it is enough to choose the parameter $a$
so that the map satisfies the Misiurewicz condition: the critical point falls
into an unstable periodic orbit after a finite number of iterates.  Misiurewicz
proved that when this condition is satisfied, the map $f$ possesses an
invariant density $\rho_0$ {\cite{misiurewicz}}.  The rest of this section
considers only Misiurewicz maps in the family (\ref{eqn:map3}).  Note that even
when the Misiurewicz condition is satisified, the action of $f$ creates a
$x^{-1/2}$ singularity at each forward image $f^n(0)$ of the critical point
$0.$ Thus the invariant density $\rho_0$ of $f$ must contain $x^{-1/2}$
singularities, as one can see in Figure {\ref{fig:map3b}}.

\begin{figure}
\begin{center}
\putgraph{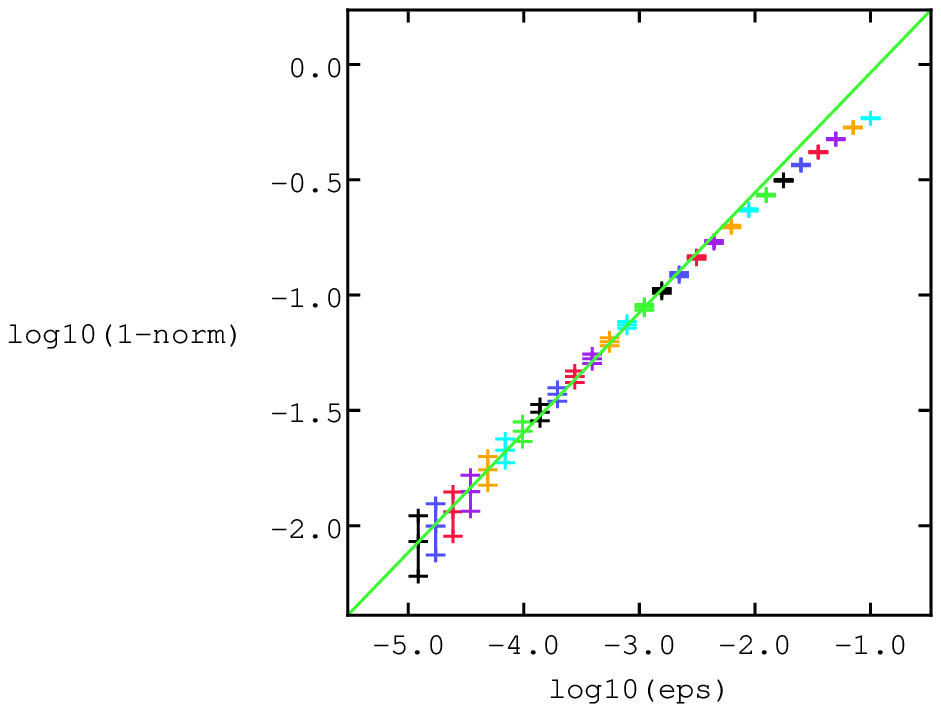}
\caption{The differences $\onenorm{\rho_\eps-\rho_0}$ and
  $\onenorm{G_\eps\rho_0-\rho_0}$ as a function of $\eps$ on a log-log graph.
  The slope of the best linear fit is {\mapthreea}; the corresponding line is
  shown as a solid line.  The slope is calculated using least-squares
  regression on the interval $\eps\in\mapthreearange.$ Note that this means
  some of the rightmost data points are discarded.  The error bars mark the
  mean square error in the computed value.}
\label{fig:map3c}
\end{center}
\end{figure}

The Perron-Frobenius operator $T_f$ is known to have a spectral gap even though
$f$ is not uniformly expanding {\cite{young3}}.  This map therefore provides a
more stringent test of the heuristic estimate than the previous two examples.
Because the invariant density $\rho_0$ contains $x^{-1/2}$ singularities, the
main effect of $G_\eps$ on $\rho_0$ is to mollify these singularities.  This is
again consistent with the picture in Figure {\ref{fig:map3b}}.  The heuristic
predicts then that
\begin{equation}
\onenorm{\rho_\eps-\rho_0}\sim\onenorm{G_\eps\rho_0-\rho_0}\sim\eps^{1/2}.
\end{equation}
The $L^1$ norms $\onenorm{\rho_\eps-\rho_0}$ were calculated for
$a=1.7152100141023...$ and plotted as a function of $\eps$ on a log-log graph
in Figure {\ref{fig:map3c}}; the error bars mark $1$ standard deviation from
the computed value.  Using least-squares regression on the interval
$\eps\in\mapthreearange$ to calculate the slope $\gamma$ yields {\mapthreea};
the corresponding line is shown as a solid line.  The data is consistent with
the prediction that $\gamma=1/2.$ More careful calculations are necessary to
decide whether $\gamma$ is exactly $1/2,$ as predicted by the heuristic.

\begin{figure}
\begin{center}
\putgraph{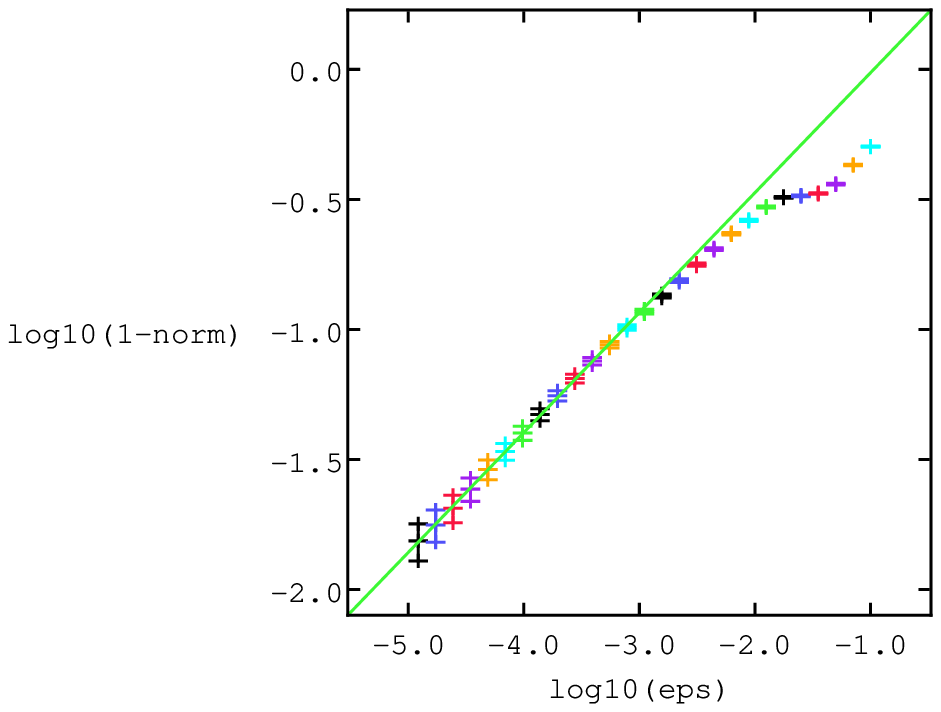}
\caption{The differences $\onenorm{\rho_\eps-\rho_0}$ and
  $\onenorm{G_\eps\rho_0-\rho_0}$ as a function of $\eps$ on a log-log graph.
  The slope is {\mapthreeb}.  The slope is calculated using least-squares
  regression on the interval $\eps\in\mapthreebrange.$ The error bars mark the
  mean square error in the computed value.}
\label{fig:map3+b}
\end{center}
\end{figure}

The computation is repeated for another Misiurewicz parameter,
$a=1.777776174649396.$ Again, using least-squares regression on the interval
$\eps\in\mapthreebrange$ to calculate the slope $\gamma$ yields {\mapthreeb}.
The results are shown in Figure {\ref{fig:map3+b}}: the data is again
consistent with the claim that $\gamma=1/2.$


\subsubsection*{Numerical method.}

Unlike the previous examples, attempts to compute invariant densities by
discretizing $T_f$ (or $T_\eps$) for (\ref{eqn:map3}) does not work
consistently: ARPACK routines will converge only for some values of $\eps$ and
not at all for others.  Sometimes ARPACK produces eigenvectors which have no
obvious connection to the dynamics or to the known form of $\rho_0.$ This
sensitive dependence on the value of $\eps$ may be related to the extreme
sensitivity of the map to the value of the parameter $a$: Misiurewicz
parameters do not form an open set; nearby values may not even have invariant
probability densities.  It is also possible that the structure of the spectral
decomposition of $T_f$ is sufficiently complex that ARPACK routines could not
produce well-resolved answers.\footnote{In contrast, Ulam's Markov chain method
should work quite well for this map.  This option was not explored because it
was not necessary.}  Fortunately, because of the $x^{-1/2}$ singularities, the
values of $\onenorm{\rho_\eps-\rho_0}$ needed here are orders of magnitude
larger than for the first two examples.  (Compare Figures {\ref{fig:map1c}},
{\ref{fig:map2c}}, {\ref{fig:map3c}}, and {\ref{fig:map3+b}}.)  Furthermore,
$T_f$ has a spectral gap, so correlations decay exponentially fast.  This
suggests that the usual time-averaging procedure, based on the ergodicity of
$f,$ can provide sufficiently accurate estimates of the probability densities
$\rho_\eps.$:
\begin{enumerate}

\item
Partition the interval $[-1,1]$ into $N$ equal-sized intervals $I_i.$

\item
Pick a random (uniform) initial condition $x_0$ and compute $x_1, x_2, ...,
x_M$ using (\ref{eqn:random-ds}).  Record the frequencies with which $x_k$
visits each of the intervals in the partition.

\item
Let $\hat{p}_i$ be the relative frequency of the $i$th interval in the
partition.  By the ergodic theorem, $\hat{p}_i$ should be approximately the
probability $\int_{I_i}{\rho_\eps(x)\ dx}.$ The $L^1$ distance between
$\rho_\eps$ and $\rho_0$ can be computed by
\begin{equation}
\label{eqn:histo-1-norm}
\onenorm{\rho_\eps-\rho_0}\approx\sum_{i=1}^{N}
{\abs{\hat{p}_i(\eps)-\hat{p}_i(0)}}.
\end{equation}
Let us denote the estimator $\sum_{i=1}^{N}
{\abs{\hat{p}_i(\eps)-\hat{p}_i(0)}}$ by $\Phi(\eps)$; it provides estimates of
$\onenorm{\rho_\eps-\rho_0}.$

\end{enumerate}
Applying the estimator $\Phi(\eps)$ to a set of noise amplitudes $\eps_i$ and
combining the results with least squares regression yields an estimator
$\hat\gamma$ of $\gamma.$ One might expect to obtain better results by adapting
the partition to better resolve the $x^{-1/2}$ singularities (their positions
can be determined {\em a priori}: they are located on the forward images of the
critical point).  However, it is easier to use uniform-sized partitions, and
they appear to be sufficient for this calculation.

The error analysis is standard.  The estimator $\Phi(\eps)$ is biased because
the finite partition used in Equation (\ref{eqn:histo-1-norm}) induces a can
only locate the zeros of $\rho_\eps-\rho_0$ up to a length scale of $N^{-1}.$
More precisely, each sign change of $\rho_\eps-\rho_0$ contributes an error of
order $N^{-2}$ to the right hand side of Equation (\ref{eqn:histo-1-norm}).  As
there are only a finite number of zero crossings, the overall bias is
$O(N^{-2}).$ The partition size $N$ in this calculation should be sufficiently
large (see below) to make the bias negligible.  The mean square error of
$\hat\gamma$ is thus assumed to be dominated by its variance.

The variance of $\hat\gamma$ can be estimated from the variance of the
$\Phi(\eps)$: the standard deviation of the estimated probability $\hat{p}_i$
is $\Delta\hat{p}_i=\sqrt{\frac{\hat{p}_i}{M}},$ where $M$ is the number of
steps taken during the course of the computation.\footnote{This assumes that
the number of points in a given interval follows a Poisson distribution.}
Summing over $i$ gives
$$\mbox{standard deviation of $\Phi(\eps)$}\leq\sqrt{\frac{N}{M}}.$$ As
$\hat\gamma$ is linear in $\log_{10}(\Phi(\eps)),$ its standard deviation is
bounded by
\begin{equation}
\frac{1}{n}\sqrt{\sum_{i=1}^{n} {\of{\frac{\log_{10}(\eps_i)-\mu}
{\sigma^2}}^2\cdot\frac{N}{M}\cdot\frac{1}{\E\brac{\Phi(\eps_i)}^2}}},
\end{equation}
where $\mu =\frac{1}{n}\sum_{i=1}^n{\log_{10}(\eps_i)},$ $\sigma^2
=\frac{1}{n}\sum_{i=1}^n{\of{\log_{10}(\eps_i)-\mu}^2},$ and
$\frac{N}{M}\cdot\frac{1}{\E\brac{\Phi(\eps)}^2}$ is an estimate of the
variance of $\log_{10}(\Phi(\eps)).$ Since the expectation values
$\E\brac{\Phi(\eps_i)}$ appear in the error bounds but are not available
exactly, they are replaced by the estimates $\Phi(\eps_i)$ themselves in the
calculations.  Note that this formula assumes that the estimates $\Phi(\eps_1)$
and $\Phi(\eps_2)$ are statistically independent if $\eps_1\neq\eps_2.$

The results shown in this section have been computed using a partition of
$N=2^{16}=65536$ intervals and $M=1.6\times 10^5\times N\approx 10^{10}$ steps.
Thus standard deviation of $\Phi(\eps)$ is on the order of $2.5\times 10^{-3}.$


\section{Intermittent maps}
\label{sec:nogap}

The examples in the previous all exhibit exponential decay of correlations.  In
this section I examine two examples whose Perron-Frobenius operators $T_f$ do
not have spectral gaps and correlation functions decay algebraically.  The
subexponential decay is caused by parts of phase space where the map $f$ is
nonexpanding.  This therefore provides a model of ``intermittency.''  See
{\cite{pomeau,liverani}}.

\subsection{Circle map with neutral fixed point}
\label{sec:map4}

\begin{figure}
\begin{center}
\putgraph{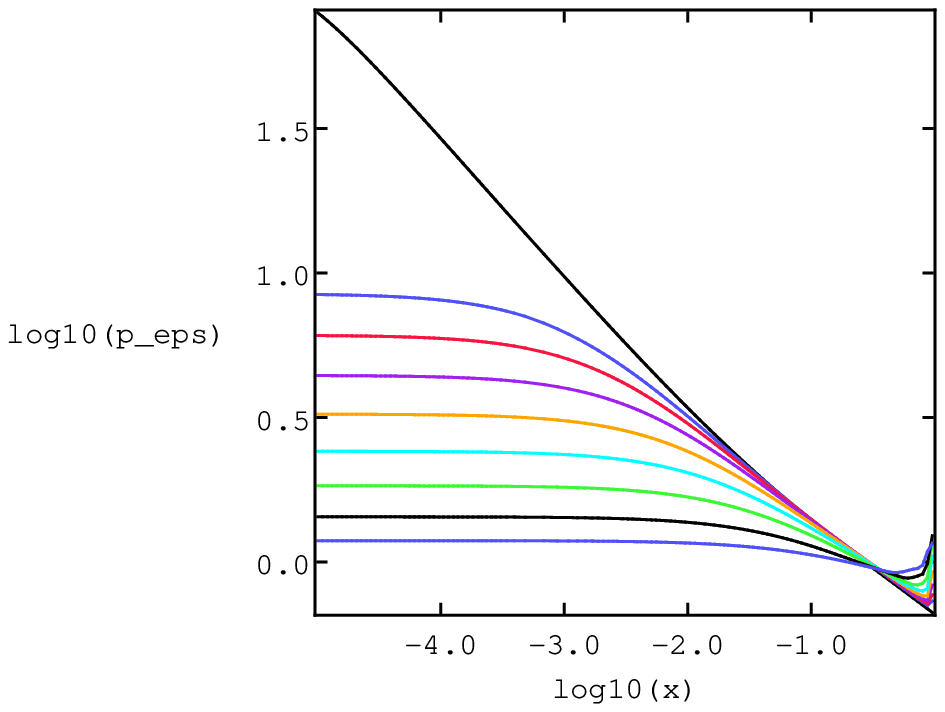}
\caption{The invariant densities $\rho_\eps$ with increasing values of $\eps.$
The map parameter is $\alpha=0.5.$ The noiseless density $\rho_0$ therefore has
a $x^{-0.5}$ singularity.  Note that this figure uses a log-log scale.}
\label{fig:map4b}
\end{center}
\end{figure}

The first of the intermittent examples is the map
\begin{equation}
\label{eqn:map4}
f(x)=\left\{\begin{array}{ll}
x+2^\alpha x^{1+\alpha}, &0\leq x\leq\frac{1}{2}\\
2x-1,                    &\frac{1}{2}<x<1\\
\end{array}\right..
\end{equation}
This is a modification of the angle-doubling map $x\mapsto 2x\mbox{ (mod $1$)}$
with a one-sided tangency to the diagonal at $x=0$: $\lim_{x\to 0+}f'(x)=1.$
Dynamically, this means that $f(x)\approx x$ when $x$ is small and positive, so
whenever $x_k$ lands near and to the right of the origin, many subsequent
iterates are required before the trajectory ``escapes'' from $0.$ The tangency
has a significant impact on the dynamics: it can be proved that correlations
decay like $n^{1-1/\alpha}$ for this map, and that the invariant density has a
$x^{-\alpha}$ singularity at $x=0$ {\cite{young4}}.  See Figure
{\ref{fig:map4b}}.

For this example, the perturbation kernel is taken to be
\begin{equation}
g(x)=\left\{\begin{array}{ll}
1,&-1\leq x\leq 0\\
0,&\mbox{ otherwise}\\
\end{array}\right..
\end{equation}
This choice is arbitrary; using (\ref{eqn:pert}) does not affect the results.

\begin{figure}
\begin{center}
\putgraph{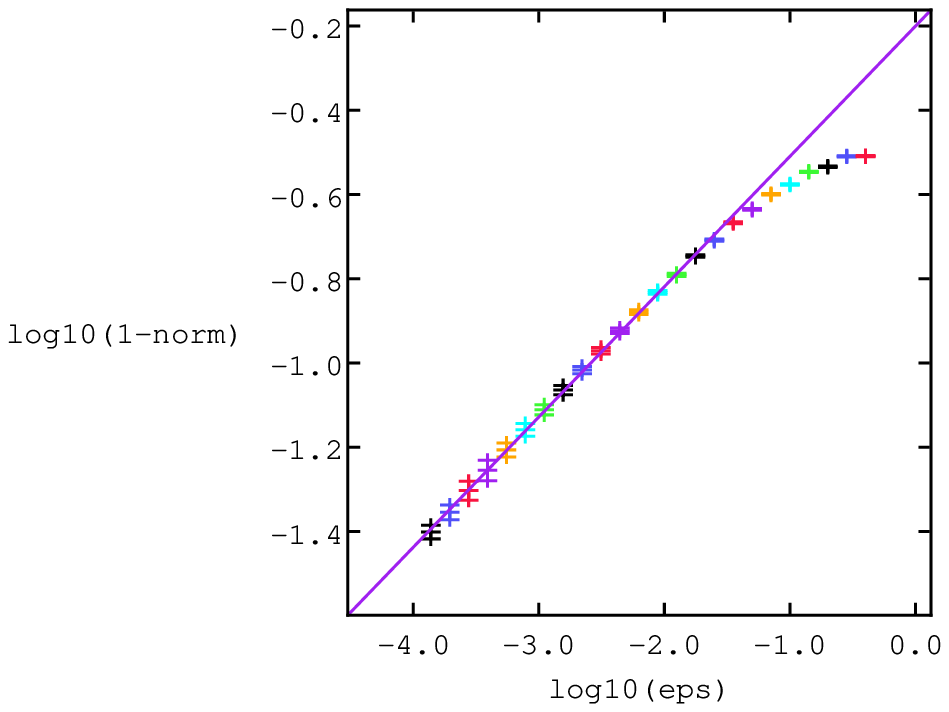}
\caption{The differences $\onenorm{\rho_\eps-\rho_0}$ and
  $\onenorm{G_\eps\rho_0-\rho_0}$ as a function of $\eps$ on a log-log graph.
  The slope is {\mapfoura} and is calculated using least-squares regression on
  the interval $\eps\in\mapfourarange.$}
\label{fig:map4c}
\end{center}
\end{figure}

As in Section {\ref{sec:map3}}, the main effect of $G_\eps$ is to smooth out
the singularity in $\rho_0.$ Thus the heuristic (\ref{eqn:heuristic1}) would
predict that
$$\onenorm{\rho_\eps-\rho_0}\sim
\onenorm{G_\eps\rho_0-\rho_0}\sim\eps^{1-\alpha}.$$ This prediction does not
change when one uses the more general form of the heuristic
(\ref{eqn:heuristic}) with $n>1.$ However, it can be seen in Figure
{\ref{fig:map4c}} that $\onenorm{\rho_\eps-\rho_0}\nsim\eps^{1-\alpha}$: in
this computation $\alpha=0.5,$ so we would expect $\eps^{0.5}$ convergence as
$\eps\to 0$.  And yet the computed exponent $\gamma$ is only {\mapfoura},
significantly smaller than $0.5.$ This result can be verified by repeating the
calculation for different values of $\alpha:$ for $\alpha=0.3$, the heuristic
predicts $\eps^{0.7}$ convergence.  The real exponent is {\mapfourb}.  And for
$\alpha=0.7$, the heuristic predicts $\eps^{0.3}$ convergence.  The real
exponent is {\mapfourc}.  The heuristic is qualitatively correct, though: as
$\alpha$ increases, $\gamma$ decreases.  Also, the computed exponents are
consistent with the fact (see Section {\ref{sec:map1}}) that the heuristic
always provides a lower bound for the error $\onenorm{\rho_\eps-\rho_0}.$


\begin{figure}
\begin{center}
\putgraph{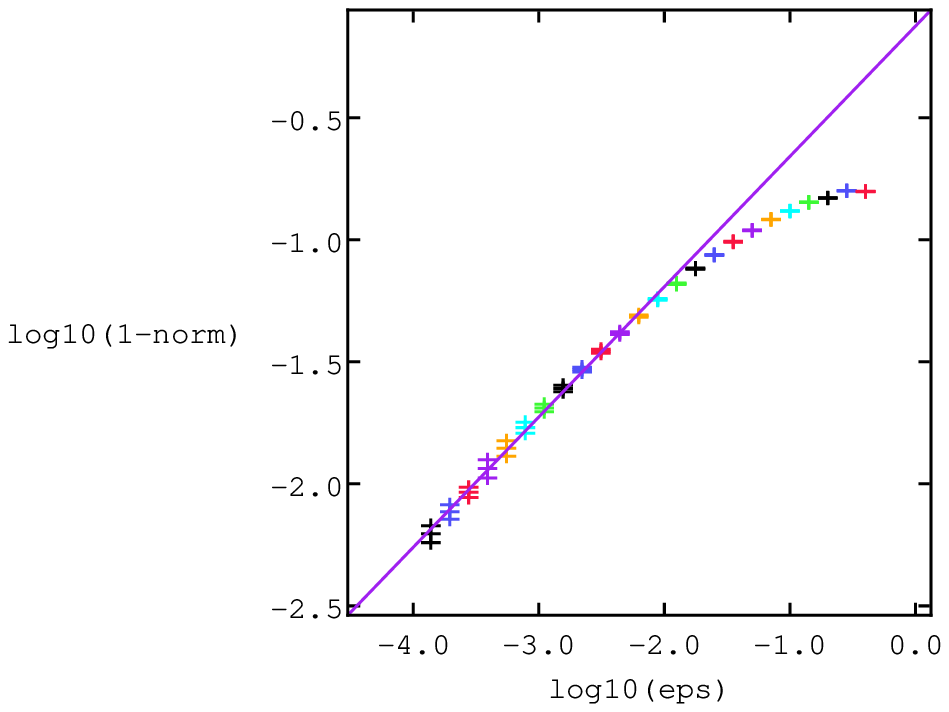}
\caption{The differences $\onenorm{\rho_\eps-\rho_0}$ and
$\onenorm{G_\eps\rho_0-\rho_0}$ as a function of $\eps$ on a log-log graph.
The slope is {\mapfourb} and is calculated using least-squares regression on
the interval $\eps\in\mapfourbrange.$}
\label{fig:map4+b}
\end{center}
\end{figure}


\begin{figure}
\begin{center}
\putgraph{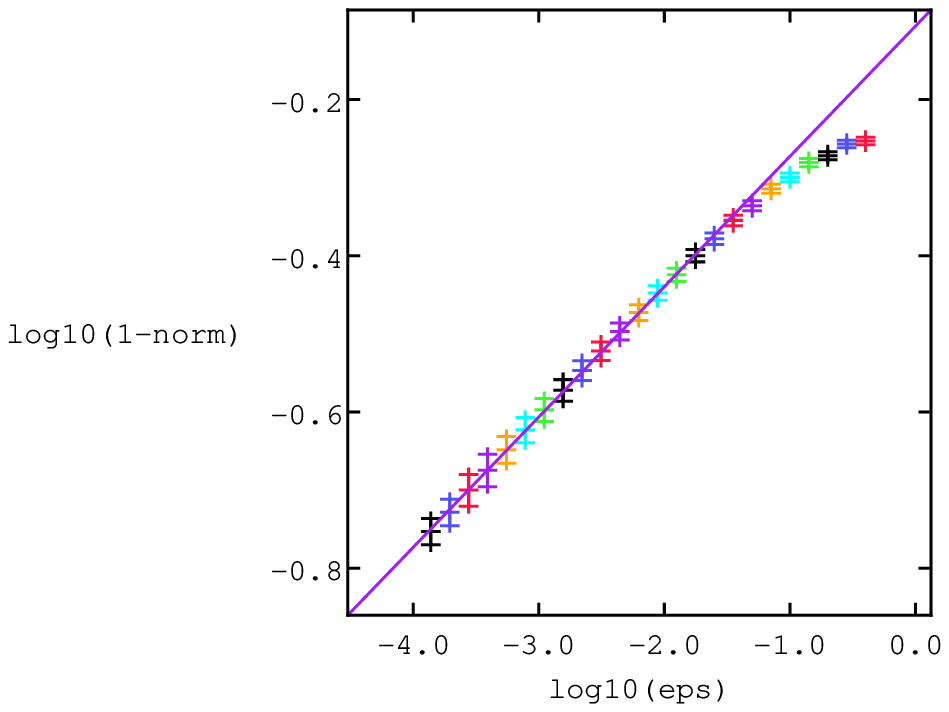}
\caption{The differences $\onenorm{\rho_\eps-\rho_0}$ and
$\onenorm{G_\eps\rho_0-\rho_0}$ as a function of $\eps$ on a log-log graph.
The slope is {\mapfourc} and is calculated using least-squares regression on
the interval $\eps\in\mapfourcrange.$}
\label{fig:map4++b}
\end{center}
\end{figure}

The numerical results indicate that the lack of a spectral gap can have a
significant, qualitative impact on the rate of convergence of $\rho_\eps$ to
$\rho_0,$ and hence on the degree of stability of $\rho_0$ under random
perturbations.

\subsubsection*{Numerical method.}

The computations for this example again rely on the finite difference scheme of
Sections {\ref{sec:map1}} and {\ref{sec:map2}}.  The only differences are:
\begin{enumerate}

\item
As in Section {\ref{sec:map2}}, it is important to avoid interpolating across
discontinuities of $f'.$ Again, this means using stencils which do not contain
$x=0$ and $x=1/2$ as interior points.  In fact, because of the singularity in
$\rho_0$ at $0,$ the origin should be excluded from the grid.

\item
In order to resolve the $x^{-\alpha}$ singularity in $\rho_0(x),$ it is
necessary to make the grid finer near $0.$ The grid I use is a hybrid between a
uniform grid and one which scales like a power law: near $x=0$ the grid
switches to one with points $\hat{x}_j\propto j^{-\alpha}, j=1,2,3,....$ This
scaling is the same as in Young's tower construction {\cite{young4}}.

\end{enumerate}
In the computations shown in this section, the polynomial interpolation of
$\hat{\rho}$ utilizes a stencil with $6$ points and the grid consists of
$N=10^4$ points.  Numerical tests show that the order of convergence is between
$5$ and $6.$ Thus the finite difference scheme is sufficiently accurate to
provide an efficient way to compute $\rho_\eps.$


\subsection{The Bunimovich stadium}
\label{sec:map5}

\begin{figure}
\begin{center}
\putgraph{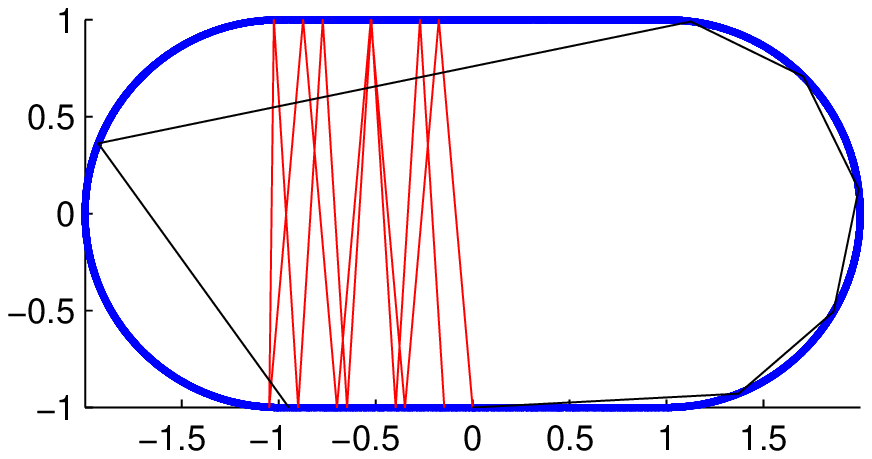}
\caption{The stadium.  Two trajectories are shown here: a perturbed ``bouncing
  ball'' orbit between the top and bottom edges, and a near ``whispering
  gallery'' orbit which wound around the right circular wing.}
\label{fig:map5a}
\end{center}
\end{figure}

\begin{figure}
\begin{center}
  \setlength{\unitlength}{2.75in}
  \begin{picture}(1.1,1.2)

  \put(0.25,1.1){\line(1,0){0.5}}
  \put(0.25,0.6){\line(1,0){0.5}}
  \put(0.25,0.85){\arc{0.5}{1.5707963267948966}{4.71238898038469}}
  \put(0.75,0.85){\arc{0.5}{-1.5707963267948966}{1.5707963267948966}}
  \put(0.76,0.86){\arc{0.5}{-1.5707963267948966}{0}}

  \put(0.95,0.85){$A$}
  \put(0.5,1.05){$B$}
  \put(0.05,0.85){$C$}
  \put(0.5,0.65){$D$}
  \put(0.5,0.85){\bf $\Omega$}

  \put(1,0.85){\circle*{0.02}}
  \put(1.02,0.85){$p_*$}
  \put(0.95,1.05){$x$}

  \put(0.75,1.1){\circle*{0.02}}
  \put(0.75,1.12){$p$}

  \put(0.75,1.1){\vector(-1,-1){0.1}}
  \put(0.69,1.07){$\theta$}

  \put(0.1,0.05){\line(1,0){0.8}}
  \put(0.1,0.55){\line(1,0){0.8}}
  \put(0.1,0.05){\line(0,1){0.5}}
  \put(0.9,0.05){\line(0,1){0.5}}

  \put(0.2,0.05){\line(0,1){0.5}}
  \put(0.4,0.05){\line(0,1){0.5}}
  \put(0.6,0.05){\line(0,1){0.5}}
  \put(0.8,0.05){\line(0,1){0.5}}

  \put(0.13,0.3){$A$}
  \put(0.85,0.3){$A$}
  \put(0.29,0.3){$B$}
  \put(0.49,0.3){$C$}
  \put(0.69,0.3){$D$}

  \put(0.49,0){$x$}
  \put(0.1,0){$0$}
  \put(0.9,0){$L$}

  \put(0.07,0.3){$\theta$}
  \put(0.07,0.05){$0$}
  \put(0.07,0.55){$\pi$}

  \end{picture}

\caption{Illustration of the coordinate system.  The angle $\theta$ ranges from
$0$ to $\pi$ while the arclength (position) $x$ ranges from 0 to $L$, the
perimeter of the stadium $\Omega.$}
\label{fig:map5c}
\end{center}
\end{figure}
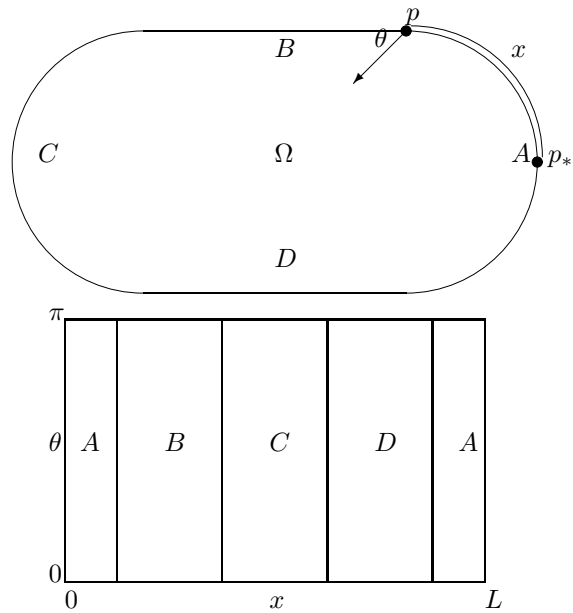

The last example is the stadium billiard: begin with a domain
$\Omega\subset\R^2$ which is the union of a rectangle and two semi-circular
ends, and consider the dynamics of a free point particle inside $\Omega$ which
collides with $\partial\Omega$ elastically (see Figure {\ref{fig:map5a}}).  The
state of the particle can therefore be represented by a point in $\Omega$
together with a unit vector.  The discussion in this section focuses on the
stadium map, which is the Poincar\'e map $f$ defined by the boundary
$\partial\Omega.$ That is, the map $f$ takes as input a pair $(x,\theta),$
where $x\in\partial\Omega$ and $\theta$ is an angle specifying the velocity of
the particle, and outputs the position and velocity of the particle {\em after}
the next collision has occurred.\footnote{This map falls outside of the
framework set up in the Introduction, but most of the general discussion there
applies to this example with only minor modifications.}  Mathematically, the
domain of $f$ is homeomorphic to $S^1\times[0,\pi],$ with the position variable
being periodic.  Note that unlike all other examples in this paper, the
billiard map is $2$-dimensional.  The map $f$ also has singularities: $Df$ is
discontinuous on the preimage of the vertical lines in Figure
{\ref{fig:map5c}}.  The expository article by Chernov and Young
{\cite{chernov}} provides a clear survey of the statistical properties of
billiards.  While they do not discuss the stadium, they do explains many
relevant ideas in terms of other billiard models.

In order to carry out numerical calculations, a specific coordinate system is
needed.  The calculations described here adhere to the following coordinates:
fix a reference point $p_*$ on $\partial\Omega$ and specify points $p$ on
$\partial\Omega$ by the length $x$ of the arc subtended by $p_*$ and $p$ in the
counterclockwise direction.  The velocity vector is specified by the angle
$\theta\in[0,\pi]$ between the vector and the tangent line to $\partial\Omega,$
again in a counterclockwise direction.  See Figure {\ref{fig:map5c}}.  In these
coordinates, the billiard map $f$ preserves an invariant density
$\rho_0(x,\theta)\propto\sin(\theta)$; the invariant density is uniform in the
$x$ (arclength) variable.  The pair $(f,\rho_0)$ is ergodic and has a positive
Lyapunov exponent {\cite{bunimovich,donnay}}.

For simplicity, rather than adding noise to both the $x$ and $\theta$
coordinates, only the $\theta$ coordinate is perturbed in this example.  The
choice of $\theta$ is natural: the geometry of the stadium suggests that small
changes in $\theta$ will be magnified very quickly; perturbations in $x$ alone
will not necessarily produce that effect.  Furthermore, perturbations in
$\theta$ destroys all metastable periodic orbits, such as the so-called
``bouncing ball'' orbits.  This ensures that the random dynamical system
(\ref{eqn:random-ds}) has a unique invariant measure.  One complexity which
arises in adding noise to $\theta$ alone is that the perturbation can no longer
be purely additive: $\theta$ must lie between $0$ and $\pi.$ In order to
satisfy this constraint, the following kernel is used:
\begin{equation}
\label{eqn:stadpert}
g_\eps(\theta_{\mbox{\tiny old}}\mapsto\theta_{\mbox{\tiny new}}) =
\left\{\begin{array}{ll}
\frac{1}{\eps+\theta_{\mbox{\tiny old}}},&0\leq\theta_{\mbox{\tiny old}}<\eps\\
\frac{1}{2\eps},&\eps\leq\theta_{\mbox{\tiny old}}\leq\pi-\eps\\
\frac{1}{\eps+\pi-\theta_{\mbox{\tiny old}}},&\pi-\eps<\theta_{\mbox{\tiny old}}\leq\pi\\
\end{array}\right..
\end{equation}
This is a somewhat arbitrary recipe; it is not clear how much of the results in
this section are due to the arbitrary nature of this specific recipe and how
much is truly intrinsic to the stadium.

\begin{figure}
\begin{center}
\putgraph{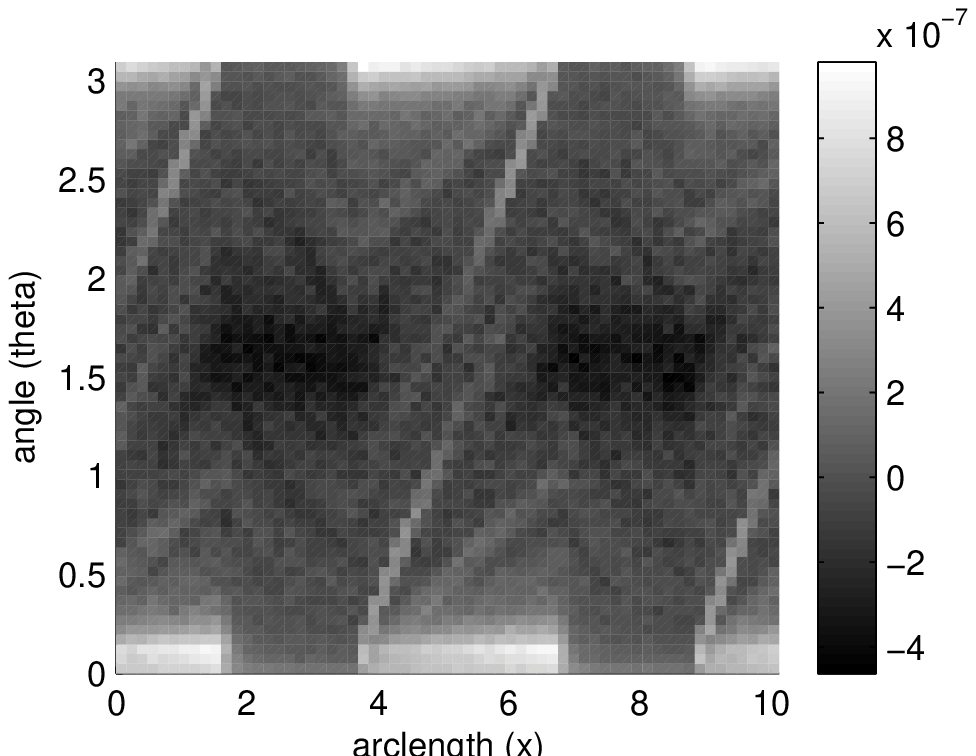}
\caption{Level sets of $\rho_\eps-\rho_0,$ with $\eps=0.1.$}
\label{fig:map5b}
\end{center}
\end{figure}

Figure {\ref{fig:map5b}} shows the difference $\rho_\eps-\rho_0.$ The
particular stadium used in the computation consists of a square $[-1,1]^2$ with
two semi-circles of radius $1$ attached.  The arclength variable therefore
ranges from $0$ up to the perimeter $L=2\pi+4$ of the stadium.  In these
coordinates, then, the vertical strips $B=[\pi/2,2+\pi/2]\times[0,\pi]$ and
$D=[3\pi/2+2,3\pi/2+4]\times[0,\pi]$ correspond to the flat edges of the
stadium, whereas the strips $A=[-\pi/2,\pi/2]\times[0,\pi]$ and
$C=[\pi/2+2,3\pi/2+2]\times[0,\pi]$ correspond to the circular ends.  As one
can see, $\rho_\eps-\rho_0$ is more negative near the middle of the strips $B$
and $D$ and more positive near the edges of the strips $A$ and $C.$ This is not
surprising: the effect of the perturbation is to decrease the amount of
probability near the vertical bouncing ball orbits, which correspond to the
parts of $B$ and $D$ near $\theta=\frac{\pi}{2},$ and the asymmetry in the
recipe (\ref{eqn:stadpert}) tends to create shallower trajectories with angles
nearer $0$ or $\pi.$

Because $\rho_0$ is smooth, the effect of applying $G_\eps$ to $\rho_0$ once is
to smooth it out further.  Thus the heuristic predicts
$\onenorm{\rho_\eps-\rho_0}\sim\eps^2.$ In this case, it is entirely possible
that the heuristic (\ref{eqn:heuristic1}) does not always apply.  Instead,
Equation (\ref{eqn:heuristic}) may be needed with $n>1.$ It is clear that, at
least for $n=2,$ the main effect of the singularties in $Df$ is to introduce
``ridges'' into $T_\eps^2\rho_0,$ {\em i.e.} lines in $S^1\times[0,L]$ along
which $T_\eps^2\rho_0$ is not differentiable in one direction.  Explicit
calculations show that such ridges contribute a term of $O(\eps^2)$ to
$\onenorm{T_\eps^2\rho_0-\rho_0}.$

\subsubsection*{Numerical method and results.}

\begin{figure}
\begin{center}
\putgraph{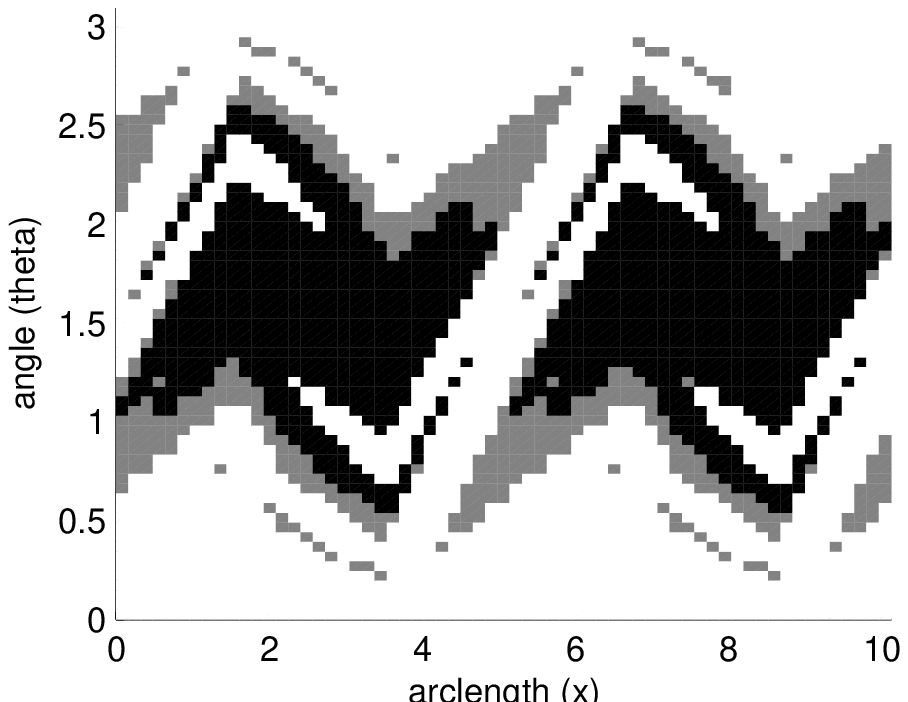}
\caption{Dynamically-generated support of the test function $\phi_\eps.$ The
  dark black set is $\set{\phi_\eps=-1},$ the white set $\set{\phi_\eps=+1}$,
  and the remainder (gray) is $\set{\phi_\eps=0}.$}
\label{fig:map5d}
\end{center}
\end{figure}

Because of the non-smoothness of the boundary $\partial\Omega,$ the stadium map
$f$ is also not smooth.  Unlike the one-dimensional case (see Sections
{\ref{sec:map2}} and {\ref{sec:map4}}), these sets of discontinuity are no
longer mere point sets but have geometric structure.  This makes the direct
discretization of $T_f$ troublesome and the finite difference method of
{\S\ref{sec:map1}}, {\S\ref{sec:map2}}, and {\S\ref{sec:map4}} difficult to
apply.  And, unlike the quadratic map in Section {\ref{sec:map3}}, the
noiseless stadium dynamics has slow decay of correlations:
$$\abs{\int_{S^1\times[0,L]} {\phi\cdot(T_f^n\psi)\cdot\rho_0}
  -\int_{S^1\times[0,L]} {\phi\rho_0}\cdot\int_{S^1\times[0,L]}
  {\psi\rho_0}}\leq cn^{-1}.$$ It implies that when $\eps$ is small but
positive, the noisy dynamics will exhibit exponential decay of correlations but
with a very large decay time constant.  This, combined with the two-dimensional
nature of the map, renders the computation of $\onenorm{\rho_\eps-\rho_0}$ by
the method of {\S\ref{sec:map3}} impractical: the error of that method is
proportional to the inverse square root of the number of samples per bin, which
would need to be quite large in this case.

Instead, notice that the basic property (\ref{eqn:onenorm}) of the $L^1$ norm
tells us that if we let
$$\phi_\eps(x) =\left\{\begin{array}{ll}
+1,&\rho_\eps-\rho_0 > 0\\
-1,&\rho_\eps-\rho_0 < 0\\
0, &\rho_\eps-\rho_0 = 0\\
\end{array}\right.$$
then
$$\onenorm{\rho_\eps-\rho_0} =\abs{\int_M{\phi_\eps\cdot\rho_\eps}
  -\int_M{\phi_\eps\cdot\rho_0}}.$$ This suggests that approximate knowledge of
the sets $E_+ =\set{\rho_\eps-\rho_0 > 0}$ and $E_- =\set{\rho_\eps-\rho_0 <
  0}$ will allow us to estimate $\onenorm{\rho_\eps-\rho_0}.$ But the geometric
structure in Figure {\ref{fig:map5b}} gives us quite a bit of information about
$E_\pm$: Let $E_+^0$ be union of the small strip
$[\pi/2,\pi/2+2]\times[0,\eps]$ and its symmetric images under the action of
the discrete symmetry group of the stadium.  The set of initial conditions
$E_+^0$ generates orbits which bounce many times with shallow angles; these are
the so-called ``whispering gallery'' orbits.  Explicit calculations show that
the brighter regions in Figure {\ref{fig:map5b}}, roughly indicating $E_+
=\set{\rho_\eps-\rho_0 > 0},$ corresponds to the union of the forward images of
$E_+^0$ under the billiard map after a few (3 or 4) iterations.  Similarly, let
$E_-^0$ denote the union of $[-\pi/2,\pi/2+2]\times[\pi/2-\eps,\pi/2+\eps]$ and
its symmetric cousins.  This set of initial conditions generates ``bouncing
ball'' orbits, and the union of the forward iterates of $E_-^0$ provides a
rough approximation of $E_- =\set{\rho_\eps-\rho_0 < 0}$ (the dark region in
Figure {\ref{fig:map5b}}).

This allows us to construct an observable $\phi_\eps$ as follows: let
$\phi_\eps$ take on the value $-1$ on the two sets near the midline of Figure
{\ref{fig:map5d}}, where $\rho_\eps-\rho_0$ is negative, and let it take on the
value $+1$ on the sets near the boundaries, where $\rho_\eps-\rho_0$ is
positive, and set $\phi_\eps=0$ elsewhere.  The quantity
\begin{equation}
\label{eqn:stadphi}
\Phi(\eps) =\abs{\iint{\phi_\eps(\rho_\eps-\rho_0)\ dx\ d\theta}}
\end{equation}
is then a lower bound of $\onenorm{\rho_\eps-\rho_0}.$ By construction,
$\phi_\eps$ should maximize $\Phi(\eps)$ so that
$\Phi(\eps)$ is as close to $\onenorm{\rho_\eps-\rho_0}$ as possible.

\begin{figure}
\begin{center}
\putgraph{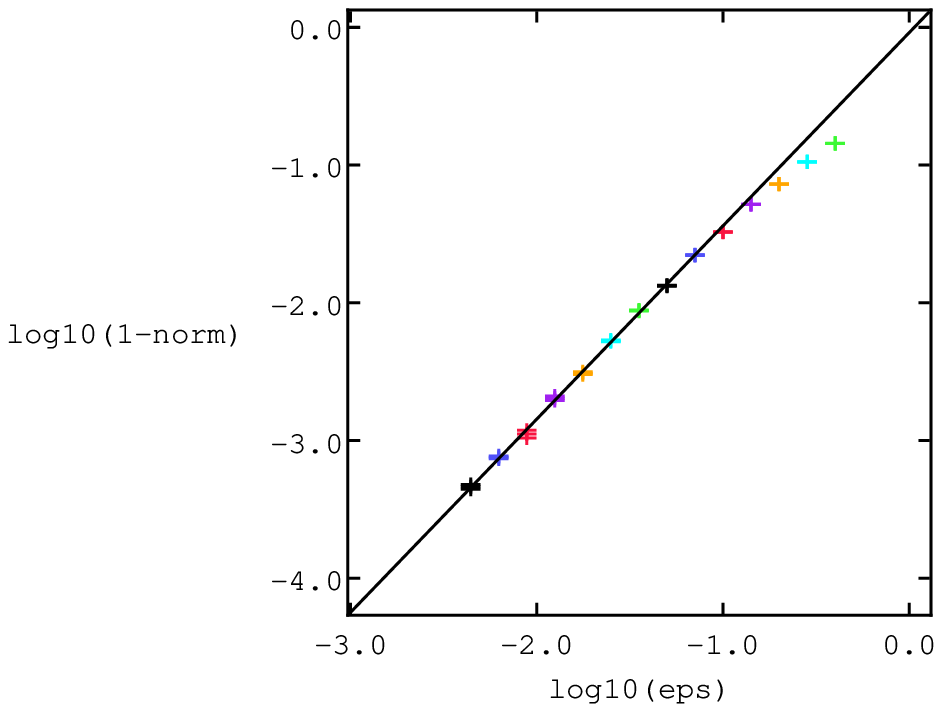}
\caption{This figure plots $\Phi(\eps)=\abs{\iint{\phi_\eps(\rho_\eps-\rho_0)\
  dx\ d\theta}}$ as a function of $\eps,$ on a log-log scale.  This provides a
  rigorous lower bound on the norm $\onenorm{\rho_\eps-\rho_0}.$ The observable
  $\phi_\eps$ is constructed to make the two quantities as close as possible.
  The slope of the best linear fit is {\mapfive} when least-squares regression
  is performed over the range $\eps\in\mapfiverange.$ This excludes points on
  the far right.  The error bars mark $1$ standard deviation.  Note that some
  of the data points have errors which are too small to be seen on this scale.}
\label{fig:map5e}
\end{center}
\end{figure}

The quantity $\Phi(\eps)$ is readily computable by averaging the values of the
observable $\phi_\eps$ over a long simulated trajectory.  The statistical error
can be estimated in a standard way {\cite{sokal}}: the standard deviation
$\Delta\Phi(\eps)$ of the estimated $\Phi(\eps)$ is bounded above by
\begin{equation}
\Delta\Phi(\eps)\leq\sqrt{\frac{2\cdot\variance(\phi_\eps)\cdot\tcorr}{N}},
\end{equation}
where the autocorrelation time $\tcorr=\sum_{n=0}^{\infty}C(n)$ and $C(n)$ is
the autocovariance function of the observable $\phi_\eps.$ The results are
shown in Figure {\ref{fig:map5e}}, where $\Phi(\eps)$ is plotted against $\eps$
on a log-log scale.  As usual, the error bars mark $1$ standard deviation.  The
slope is approximately {\mapfive} when fitted over the range
$-2.5\leq\log_{10}(\eps)\leq-1.$ As $\Phi(\eps)$ only provides a lower bound on
$\onenorm{\rho_\eps-\rho_0},$ the main conclusion of this computation is that
$\onenorm{\rho_\eps-\rho_0}$ cannot decay faster than $\eps^{1.4}.$




For Figure {\ref{fig:map5e}}, $\Phi(\eps)$ is computed, when $\eps>0.00625,$
using a trajectory consisting of $5.1\times 10^8$ steps.  When $\eps\leq
0.00625,$ the computation uses $8.1\times 10^9.$ This is necessary because in
order to determine the exponent $\gamma,$ one must compute
$\log_{10}(\Phi(\eps))$ accurately.  But the absolute error in
$\log_{10}(\Phi(\eps))$ is proportional to the relative error in $\Phi(\eps),$
so it is necessary to use more steps when $\Phi(\eps)$ is small.


\section{Concluding Remarks}

The calculations described in this paper lead to many intriguing questions.  In
addition to the question of how one might formulate and prove a precise version
of the heuristic estimate (\ref{eqn:heuristic}), there is also the question of
how to correctly predict the scaling of $\onenorm{\rho_\eps-\rho_0}$ in simple
intermittent systems.  For this question, large deviations theory
{\cite{freidlin,reimann}} or renormalization techniques {\cite{hirsch}} may be
relevant.

Invariant densities represent the most regular type of invariant measures.  In
most dissiptive systems, invariant measures are supported on attractors of zero
measure.  Among such singular invariant measures, the best-behaved are SRB
measures: they represent the ``nicest'' invariant measures one can hope to have
in dissipative chaotic dynamics.  It is natural ask the corresponding question
of convergence rates for SRB measures, say the rate at which $\mu_\eps$
converges to $\musrb$ in the total variation norm.  The answer, however, is not
so apparent.  It seems sensible to conjecture that, at least in uniformly
hyperbolic systems, the rate of convergence may be determined by the regularity
of the conditional densities of $\musrb$ along unstable manifolds.

Another set of open questions have to do with dimension.  All the examples
considered in this paper exist in low-dimensional spaces.  Are there other
factors which can affect the rate of convergence in higher dimensions which
cannot be seen in low dimensions?

The numerical calculations described in this paper employed a variety of
methods.  The convergence properties of these numerical methods is not yet
completely understood and await deeper analysis.  In particular, the extent to
which a discrete transfer operator $\hat{T}$ captures the detailed spectral
structure of $T_f$ and the effect that this has on the convergence of numerical
eigenproblem routines is far from understood, as indicated by the sensitivity
of the computed density for the Misiurewicz maps of {\S\ref{sec:map2}} to
changes in the parameters $a$ and $\eps.$

Finally, as mentioned earlier, the results described here may be relevant for
numerical studies of dynamical systems with intermittent, metastable behavior.
Noise can help reduce initialization bias in long-time numerical simulations
while introducing controllable errors.  This may be particularly useful in
intermittent systems of the type examined in Sections {\ref{sec:map4}} and
{\ref{sec:map5}}.  A first step in this direction was made in {\cite{lin1}},
but a systematic exploration is required to understand these ideas.  In
particular, to carry out this idea in practice will require algorithms which
can cope with the additional complexity of separatrices and multiple ergodic
components.


\section{Acknowledgements}

It is a pleasure to thank Lai-Sang Young for her generous help with this
project.  I am also grateful to Alex Barnett, Toufic Suidan, and George
Zaslavsky for helping to improve the exposition and for many helpful and
pleasant conversations.  This work is supported by the National Science
Foundation through a Mathematical Sciences Research Postdoctoral Fellowship.

\bibliographystyle{siam}
\bibliography{pert1}

\begin{thebibliography}{10}

\bibitem{baladi}
{\sc V.~Baladi}, {\em Positive Transfer Operators and Decay of Correlations},
  World Scientific, 2000.

\bibitem{baladi1}
{\sc V.~Baladi and L.-S. Young}, {\em On the spectra of randomly perturbed
  expanding maps}, Communications in Mathematical Physics, 156 (1993).

\bibitem{blank1}
{\sc M.~Blank and G.~Keller}, {\em Random perturbations of chaotic dynamical
  systems: stability of the spectrum}, Nonlinearity, 11 (1998), pp.~1351--1364.

\bibitem{bunimovich}
{\sc L.~A. Bunimovi{\v{c}}}, {\em The ergodic properties of certain billiards},
  Funkcional. Anal. i Prilo\v zen., 8 (1974), pp.~73--74.

\bibitem{chernov}
{\sc N.~Chernov and L.~S. Young}, {\em Decay of correlations for {L}orentz
  gases and hard balls}, in Hard ball systems and the Lorentz gas, vol.~101 of
  Encyclopaedia Math. Sci., Springer, Berlin, 2000, pp.~89--120.

\bibitem{donnay}
{\sc V.~J. Donnay}, {\em Using integrability to produce chaos: billiards with
  positive entropy}, Communications in Mathematical Physics, 141 (1991).

\bibitem{eckmann}
{\sc J.-P. Eckmann and D.~Ruelle}, {\em Ergodic theory of chaos and strange
  attractors}, Reviews of Modern Physics, 57 (1985), pp.~617--656.

\bibitem{freidlin}
{\sc M.~I. Freidlin and A.~D. Wentzell}, {\em Random Perturbations of Dynamical
  Systems}, Springer-Verlag, 1998.

\bibitem{hirsch}
{\sc J.~E. Hirsch, M.~Nauenberg, and D.~J. Scalapino}, {\em Intermittency in
  the presence of noise: a renormalization group formulation}, Physics Letters,
  87A (1982).

\bibitem{keane}
{\sc M.~Keane, R.~Murray, and L.-S. Young}, {\em Computing invariant measures
  for expanding circle maps}, Nonlinearity, 11 (1998), pp.~27--46.

\bibitem{kifer}
{\sc Y.~I. Kifer}, {\em Small random perturbations of certain smooth dynamical
  systems}, Izv. Akad. Nauk SSSR Ser. Mat., 38 (1974), pp.~1091--1115.

\bibitem{lasota}
{\sc A.~Lasota and J.~A. Yorke}, {\em On the existence of invariant measures
  for piecewise monotonic transformations}, Transactions of the American
  Mathematical Society, 186 (1973), pp.~481--488.

\bibitem{arpack}
{\sc R.~B. Lehoucq, D.~C. Sorensen, and C.~Yang}, {\em ARPACK User's Guide:
  Solution of Large-Scale Eigenvalue Problems with Implicitly Restarted Arnoldi
  Methods}, SIAM, 1998.

\bibitem{lin1}
{\sc K.~K. Lin}, {\em Random perturbations of {SRB} measures and numerical
  studies of chaotic dynamics}, PhD thesis, University of California at
  Berkeley, 2003.
\newblock Lawrence Berkeley National Lab Technical Report 53522.

\bibitem{liverani}
{\sc C.~Liverani, B.~Saussol, and S.~Vaienti}, {\em A probabilistic approach to
  intermittency}, Ergodic Theory \& Dynamical Systems, 19 (1999), pp.~671--685.

\bibitem{misiurewicz}
{\sc M.~Misiurewciz}, {\em Absolutely continuous measures for certain maps of
  an interval}, Publications Math\'ematiques de l'IH\'ES, 53 (1981),
  pp.~17--51.

\bibitem{pomeau}
{\sc Y.~Pomeau and P.~Manneville}, {\em Intermittent transition to turbulence
  in dissipative dynamical systems}, Communications in Mathematical Physics, 74
  (1980), pp.~189--197.

\bibitem{reimann}
{\sc P.~Reimann and P.~Talkner}, {\em Invariant densities for noisy maps},
  Physical Review A, 44 (1991), pp.~6348--6363.

\bibitem{sokal}
{\sc A.~D. Sokal}, {\em {Monte Carlo} methods in statistical mechanics:
  Foundations and new algorithms}, in Functional Integration (Carg\`ese, 1996),
  vol.~361 of NATO Adv. Sci. Inst. Ser. B Phys., Plenum, 1997, pp.~131--192.

\bibitem{young3}
{\sc L.-S. Young}, {\em Decay of correlations for certain quadratic maps},
  Communications in Mathematical Physics, 146 (1992), pp.~123--138.

\bibitem{young}
\leavevmode\vrule height 2pt depth -1.6pt width 23pt, {\em Ergodic theory of
  differentiable dynamical systems}, in Real and Complex Dynamics, vol.~464 of
  NATO Adv. Sci. Inst. Ser. C Math. Phys. Sci. (Hiller\o d, 1993), Kluwer
  Academic Publishers, 1995, pp.~293--336.

\bibitem{young4}
\leavevmode\vrule height 2pt depth -1.6pt width 23pt, {\em Recurrence times and
  rates of mixing}, Israel Journal of Mathematics,  (1999).

\end{thebibliography}
\end{document}